\documentclass[twoside,draft]{article}
 \usepackage{t1enc}
\usepackage{amsfonts,amssymb}
 \usepackage{latexsym,amsmath,amscd}
 \usepackage{enumerate}
\usepackage{color}
\numberwithin{equation}{section}
 \newcommand{\be}{\begin{equation}}
\newcommand{\ee}{\end{equation}}
\newcommand{\bde}{\begin{displaymath}}
\newcommand{\ede}{\end{displaymath}}
\newcommand{\beq}{\begin{eqnarray*}}
\newcommand{\eeq}{\end{eqnarray*}}
\newcommand{\beqa}{\begin{eqnarray}}
\newcommand{\eeqa}{\end{eqnarray}}

\newtheorem{Theorem}{Theorem}[section]
\newtheorem{Proposition}[Theorem]{Proposition}
\newtheorem {Cor}[Theorem]{Corollary}
\newtheorem {pro}[Theorem]{Proposition}

\newtheorem {Lemma}[Theorem]{Lemma}
\newtheorem {rem}[Theorem]{Remark}
\newtheorem {rems}[Theorem]{Remarks}
\newtheorem {com}[Theorem]{Comment}
\newtheorem {coms}[Theorem]{Comments}
\newtheorem {Definition}[Theorem]{Definition}

\newcommand{\bnota}{\begin{notation} \rm } \newcommand{\enota}{\end{notation}}
\newcommand{\bcom}{\begin{com} \rm } \newcommand{\ecom}{\end{com}}
\newcommand{\bcoms}{\begin{coms} \rm } \newcommand{\ecoms}{\end{coms}}
\newcommand {\bdef}{\begin{Definition}}
\newcommand {\edefi}{\end{Definition}}
\newcommand {\bl}{\begin{Lemma}}
\newcommand {\el}{\end{Lemma}}
\newcommand {\bethe}{\begin{Theorem}}
\newcommand {\eethe}{\end{Theorem}}
\newcommand {\bp}{\begin{pro}}
\newcommand {\ep}{\end{pro}}
\newcommand {\bcor}{\begin{Cor}}
\newcommand {\ecor}{\end{Cor}}
 \newcommand {\brem }{\begin{rem} \rm }
\newcommand {\erem }{\end{rem}}
 \newcommand {\brems }{\begin{rems} \rm }
\newcommand {\erems }{\end{rems}}

\def \ind{1\!\!1}

\def\cA{{\cal A}}
\def\G{{\cal G}}

\def\cH{{\cal H}}
\def\M{{\cal M}}
\def\N{{\cal N}}

\def\ff{{\mathbb F}}
 \def\F{{\cal F}}

\def\gg{{\mathbb G}}
\def\kk{{\mathbb K}}
\def\Q{{\mathbb Q}}
\def\P{{\mathbb P}}
\def\R{{\mathbb R}}
\def\E {\mathbb{E}}

\def\graph#1{[\! [ #1]\! ]}

\def \id{1\!\!1}
\def \Rbrack {[\![}
\def \Lbrack {]\!]}

\newcommand{\is}{\centerdot}
\def\cro#1{\langle #1\rangle}
\def \wh  {\widehat }
\def \wt  {\widetilde }
\def \proof {{\sc{Proof:}}~}
\def \finproof {\hfill $\square$  \\ }
 \def \mart {\stackrel{{\rm{ \gg -mart}}}{=}}

\title{Arbitrages    in a Progressive Enlargement Setting}

\author{ {Anna
Aksamit\footnote{Laboratoire Analyse et Probabilit\'es,
Universit\'e d'Evry Val d'Essonne,  Evry, France}, Tahir
Choulli\footnote{Mathematical and Statistical Sciences Depart.,
University of Alberta, Edmonton, Canada},    Jun Deng
\,\footnotemark[2] and Monique Jeanblanc\,\footnotemark[1]
}}

\date{\today}

\begin{document}
\maketitle

\begin{abstract}
This paper completes the analysis of Choulli et al. \cite{acdj} and  contains two principal contributions. The first contribution consists in providing and analysing many practical examples of market models that admit classical arbitrages while they preserve the  No Unbounded Profit with Bounded Risk (NUPBR hereafter) under random horizon and when an honest time is incorporated for particular cases of models. For these markets, we calculate explicitly the arbitrage opportunities. The second
contribution lies in providing simple proofs for the stability of  the No Unbounded Profit with Bounded Risk under random horizon and after honest time satisfying additional important condition for particular cases of
models.
\end{abstract}

\section{Introduction}

This paper studies  a financial market in which some assets, with
prices adapted with respect to a reference filtration $\ff$, are traded. One
then assumes that an agent has some extra information, and may use
strategies adapted to a larger filtration $\gg$. This  extra
information is modeled by the knowledge of some random time
$\tau$, when this time occurs. We restrict our study to progressive enlargement of filtration setting,  and we pay a particular attention to honest times.
 Our goal is to detect if the
knowledge of $\tau$ allows for some arbitrage, i.e., if  using
$\gg$-adapted strategies, the agent  can make profit.

\noindent  In this paper we consider two main notions of no-arbitrage, namely no classical arbitrage and  No Unbounded Profit with Bounded Risk. To the best
of our knowledge, there are no references for the case of
classical arbitrages in a general setting. The goal of the present
paper is firstly to introduce the problem, to solve it in some
specific cases and to give some explicit examples of  classical
arbitrages (with a proof different   from the one in \cite{fjs}),
and secondly to give,  in some specific models, an easy proof
of  No Unbounded Profit with Bounded Risk condition.

\noindent  In the
case of honest times avoiding stopping times in a continuous  filtration,  the same problem was studied in Fontana et al. \cite{fjs}  where the authors have investigated several kinds of arbitrages. We refer the reader to that paper for an extensive list of related results in the literature.

\noindent The paper is organized as follows: Section 2 presents the problem and recalls some definitions and results on  arbitrages and progressive enlargement of filtration. In Section 3 we study two classical situations in enlargment of filtration theory, namely immersion and positive density hypothesis cases.
Section 4 concerns honest times, and we show that, in case of a complete market, there exist classical arbitrages before and after  the honest time, and we give a way to construct these arbitrages.  This  fact is illustrated by many examples, where we exhibit  these arbitrages in a closed form. In Section 5, we study some examples of non-honest times. In Section 6, we study NUPBR condition   before a random
time and after an honest time, in some specific examples.

\section{General framework}

We consider a filtered probability space $(\Omega, \cA, \ff,\P)$
where the filtration $\ff$ satisfies the usual hypotheses and
$\F_\infty \subset \cA$, and a random time $\tau$ (i.e., a
positive  $\cA$-measurable random variable). We assume that the
financial market where  a risky asset with price  $S$ (an
$\ff$-adapted positive process) and a riskless asset $S^0$
(assumed, for simplicity, to have a  constant price so that the
risk-free interest rate is null) are traded is arbitrage free.
More precisely, without loss of generality we assume that $S$ is a $(\P,\ff)$-(local)
martingale.  In this paper, the
horizon is equal to $\infty$.

\noindent We denote by $\gg$ the progressively enlarged filtration of $\ff$
by $\tau$, i.e., the smallest right-continuous filtration
that contains $\ff$ and makes $\tau$ a stopping time defined as
$$\G_t= \cap _{\epsilon
>0} \F_{t+\epsilon}\vee \sigma (\tau \wedge (t+\epsilon)).$$
We recall that  $(\cH^\prime)$ hypothesis is said to hold between two filtrations $\ff$ and $\gg$ where $\ff \subset \gg$ if any $\ff$-martingale is a $\gg$-semimartingale.
For a semimartingale $X$ and a predictable process $H$,  we use the notation $H\centerdot  X$ for the stochastic integral $\int_0^\cdot H_sdX_s$ { when it exists}.

\noindent We start by an   elementary remark:  assume that there are no
arbitrages  using  $\gg$-predictable strategies and that $\P$ is the unique probability measure making
{$S$} an $\ff$-martingale. So, in particular, the
$(S, \ff)$ market is complete (i.e., the market where $(S,S^0)$ are
traded). Then, roughly speaking,  $S$ would be
a $(\Q,\gg)$-martingale for some equivalent martingale measure $\Q$,  hence  would be
also a $(\Q,\ff)$-martingale\footnote{Note that if $S$ is a
$(\Q,\gg)$-strict local martingale for some equivalent martingale measure $\Q$,  one
can not deduce that it is also a $(\Q,\ff)$-local martingale} and
$\Q$ will coincide with $\P$ on $\ff$.
 This implies that any $(\ff,\Q)$-martingale is a $(\gg,\Q)$-martingale.

\noindent Another trivial remark is that, in the particular case where $\tau$ is an $\ff$-stopping time, the enlarged filtration and the reference filtration are the same. Therefore, no-arbitrage conditions hold before and after $\tau$.

\subsection {Illustrative examples} \label{illus}

We study here two basic examples, in order to show in a first step
how arbitrages can occur in a Brownian filtration, and in a second
step that discontinuous models present some difficulties.

\subsubsection{Brownian case}

Let $dS_t=S_t\sigma dW_t$, where $W$ is a Brownian
motion and $\sigma$ a constant, be the price of the risky asset.
This martingale $S$ goes to $0$ a.s. when $t$ goes to infinity, hence
the random time $\tau = \sup\{ t\,:\, S_t =S^* \}$  where
$S^*=\sup_{s\geq 0 } S_s$ is a  finite honest time, and obviously
leads to an arbitrage before $\tau$: at time $0$, buy one share of
$S$ (at price $ S_0$), borrow  $ S_0$, then, at time $\tau$,
reimburse the loan $ S_0$ and sell the share of the asset at price
$S_\tau$. The gain is $ S_\tau - S_0>0$ with an initial wealth null.
There are also arbitrages after $\tau$: at time $\tau$, take a short position on $S$, i.e., hold a self financing portfolio with value $V$ such that $dV_t=-dS_t, V_\tau = 0$. Usually shortselling positions are not admissible, since $V_t=- S_t+S_\tau $ is not bounded below. Here $- S_t+S_\tau$ is positive, hence shortselling is an arbitrage opportunity.

\subsubsection{Poisson case}

Let $N$ be a  Poisson process   with intensity $\lambda$   and $M$ be its   compensated
martingale.  We define  the price process $S$ as  $dS_t=S_{t-}
\psi dM_t, S_0=1$ with $\psi$ is a constant satisfying $\psi
>-1$ and $\psi \neq 0$, so that
$$S_t=  \exp (-\lambda \psi t+ \ln (1+\psi)N_t)\,.$$
Since $\frac{N_t}{t}$ goes to
$\lambda$ a.s. when $t$ goes to infinity and $\ln(1+\psi)-\psi <0$,
$S_t$ goes to $0$ a.s. when $t$ goes to infinity. The random time
$$\tau = \sup\left \{ t\,:\, S_t =S^*\right \}$$ with
$S^*=\sup_{s\geq 0 } S_s$ is a finite honest time.

\noindent If $\psi>0$,  then $S_\tau \geq S_0$ and an  arbitrage opportunity
is realized at time $\tau$, with a long position in the stock.
If $\psi <0$, then the arbitrage is not so obvious.
We shall discuss that with more details in Section \ref{poisson}.

\noindent  There are arbitrages
 after $\tau$, selling  at time $\tau$ a contingent claim with payoff $1$, paid at
 the first time $\vartheta$ after $\tau$ when  $S_t> \sup_{s\leq \tau}S_{s}$. For $\psi>0$, it reduces to    $S_\tau=\sup_{s\leq \tau}S_s$, and,   for $\psi<0$, one has    $S_{\tau-}=\sup_{s\leq \tau}S_s$. At time $t_0=\tau$, the non informed buyer
 will agree to pay a positive price, the informed seller knows
 that the exercise will be never done.

\subsection{Admissible portfolio and arbitrages
opportunities} \label{defarbitr}

 In this section, we recall the basic definitions on
arbitrages, and we give sufficient conditions for  no arbitrages  {in a market with zero interest rate}. We refer to  \cite{fjs} for details.

\noindent Let
 $  \kk $ be one of the filtrations $ \bigl\{ \ff ,  \gg \bigr\}$.
  Note that, in order that the  integral $\theta \centerdot S$  has a meaning for  a $\gg$ predictable process $\theta$, one needs that $S$ is  a $\gg$-semimartingale. This  requires (on $\{t>\tau\}$) some hypotheses on $\tau$.

\noindent For $a\in\R_+$, an element
$\theta\in L^{\kk}\left(S\right)$ is said to be an
 $a$-admissible
$\kk$-strategy  if $\left(\theta\centerdot S\right)_{\infty}
:=\lim_{t\rightarrow\infty}\left(\theta\centerdot S\right)_t$
exists and $V_t(0,\theta):=\left(\theta\centerdot
S\right)_t\geq-a$ $\P$-a.s. for all $t\geq 0$. We denote by
$\cA^{\kk}_a$ the set of all $a$-admissible $\kk$-strategies. A
process $\theta\in L^{\kk}\left(S\right)$ is called an
\emph{admissible $\kk$-strategy} if
$\theta\in\cA^{\kk}:=\bigcup_{a\in\R_+}\!\cA^{\kk}_a$.

\noindent An admissible strategy yields an
 Arbitrage
Opportunity  if $V\left(0,\theta\right)_{\infty}\geq 0$ $\P$-a.s.
and $\P\bigl(V\left(0,\theta\right)_{\infty}>0\bigr)>0$. In order
to avoid confusions, we shall call these arbitrages \emph{classical
 arbitrages}.
If there exists no such $\theta\in\cA^{\kk}$ we say that the
financial market $\M(\kk):=\left(\Omega, \kk,\P;S\right)$
satisfies the  No Arbitrage (NA) condition.

\noindent No Free Lunch with Vanishing Risk  (NFLVR) holds in the financial
market $\M(\kk)$ if and only if there exists an equivalent
martingale measure in $\kk$, i.e., a probability measure $\Q$, such that  $\Q\sim \P$  and  the process
$S $ is a $(\Q,\kk)$-local martingale. If   NFLVR holds,   there are no classical
arbitrages.

\noindent A non-negative
$\mathcal{K}_\infty$-measurable random variable $\xi$ with
$\P\left(\xi>0\right)>0$ yields an
 {Unbounded Profit with Bounded Risk  if for all $x>0$ there exists an element
$\theta^x\in\cA^{\kk}_x$ such that
$V\left(x,\theta^x\right)_{\infty}:= x+(\theta^x \centerdot
S)_\infty\geq\xi$ $\P$-a.s. If there exists no such random
variable, we say that  the financial market $\M(\kk)$ satisfies
the No Unbounded Profit with Bounded Risk (NUPBR) condition.

\noindent We
recall that {NFLVR} holds if and only if both NA and NUPBR hold (see \cite{ds} Corollary 3.4 and Proposition 3.6, \cite{KK}).

\noindent A strictly positive $\kk$-local
martingale $L=\left(L_t\right)_{t\geq 0}$ with $L_0=1$ and
$L_{\infty}>0$ $\P$-a.s. is said to be a
 local martingale
deflator in $(S, \kk)$ on the time horizon   $\left[0,\varrho\right]$
if the process $LS^{\varrho}$ is a  $\kk$-local martingale; here
$\varrho$ is a  $\kk$-stopping time.
The important result giving the characterisation of NUPBR condition for strictly positive price process is stated in Theorem 4.12 in \cite{KK} and then generalized in Theorem 5 in \cite{ta}. We recall it here.
\bethe
\label{deflator}
Let $S$ be a strictly positive $\kk$-semimartingale. Then,
the NUPBR condition holds in $\kk$ if and only if
there exists a local martingale deflator in $\kk$ .
\eethe

\subsection{Enlargement of filtration results}
We now recall some basic results on progressive enlargement of
filtrations. The reader can refer to  Jeulin \cite{jeu} and Jeulin and Yor  \cite{jy:gf}  for more
information.

\noindent Let $\tau$ be a   random time, i.e., a  positive  random variable.
We define the right-continuous with left limits
$\ff$-supermartingale
$$
Z_t := \P\left(\tau >t\ \big|\ {\cal F}_t\right).
$$
Note that $Z_0=1$ if $\P(\tau>0)=1$.
  The optional decomposition of $Z$ leads to
an  important $\ff$-martingale that we denote by $m$,  given by
\be\label{m} m := Z+A^{o },\ee
 where $A^{o }$ is the $\mathbb F$-dual optional projection\footnote{See Appendix for the definition  if needed} of $A:=\ind_{\Rbrack\tau,
 \infty\Rbrack}$ (so $A^o$ is a non-decreasing   process). Note that $m$ is non-negative: indeed $m_t=\E(A^o _\infty + Z_\infty \vert \F_t)$.

\noindent A second important $\ff$-supermartingale, defined
through
$$  \widetilde Z_t := \P\left(\tau \geq t\ \big|\ {\cal
F}_t\right),$$ will play a particular r\^ole in the following. One
has $\tilde Z= Z+ \Delta A^{o }$, hence the supermartingale $\wt
Z$   admits a decomposition as
 \be\label{mt}   \wt Z= m- A^{o }_-\,.\ee

\noindent We start with the following obvious (but useful) result
\bl  \label{lemmem} Assume that  the  financial market $(S, \ff)$
is complete and let $\varphi$ be the $\ff$-predictable process
satisfying $m=1+\varphi\centerdot S$.
 If $m_\tau\geq 1 $   and $\P(m_\tau> 1)>0$, then, the $\gg$-predictable process $\varphi\id_{\Rbrack0, \tau\Lbrack}$ is a classical  arbitrage strategy in the market ''before $\tau$'', i.e., in $(S^\tau, \gg)$.\el
\proof The $\ff$-predictable process $\varphi$ exists due to the
market completeness. Hence $ \ind_{ \Rbrack 0,\tau\Lbrack}\varphi$
is a $\gg$-predictable  admissible self-financing strategy with
initial value 1 and  final value $m_\tau-1$ satisfying
$m_\tau-1\geq 0$ a.s. and $\P(m_\tau-1>0)>0$, so it is a classical
arbitrage strategy in $(S^\tau, \gg)$.\finproof

 \subsubsection{Decomposition formula  before $\tau$}
In a first step, we   restrict our attention to what happens before
$\tau$. Therefore,  we do not require any extra hypothesis on
$\tau$, since, for any random time $\tau$, any  $\ff$-martingale stopped at $ \tau$ is a $\gg$-semimartingale, as established by Jeulin \cite[Prop. (4,16)]{jeu}: to any
${\mathbb  F }$-local martingale $X$, we associate the ${\mathbb G
}$-local  martingale $\wh X$ (stopped at time $\tau$)
\begin{equation}
\label{jeulinavant} \wh X_t  := X_t^\tau - \int_{0}^{t\wedge \tau} {
\frac{d \langle X, m\rangle^{\mathbb F}_s }{Z_{s-}} } ,
\end{equation} where, as usual, $X^\tau$ is the stopped process defined as $X^\tau_t=X_{t\wedge \tau}$.

\noindent An interesting case  is the one of pseudo-stopping times.    We
recall that a random time $\tau$ is a pseudo-stopping time   is
any $\ff$-martingale stopped at $\tau$ is a $\gg$-martingale (see
\cite{ny}).  This is equivalent to the fact that the
$\ff$-martingale $m$ is constantly equal to 1.

\subsubsection{Honest times and decomposition formula  after $\tau$}
 We need to impose conditions on $\tau$ such that
the ($\ff$-martingale) price process $S$ is a
$\gg$-semimartingale, so that  one can define stochastic
integrals  of $\gg$ predictable processes with respect to $S$. In
this paper, we are not interested by necessary and sufficient
conditions, these ones being far from tractable (see \cite[III,
2,c]{jeu}). Instead we focus here on honest times.

\noindent    \bethe \label{honest} Let $\tau$ be a random time.
Then, the following conditions are equivalent:\\
(i) The random time $\tau$ is honest, i.e., for each $t\geq 0$,
there exists an $\F_t$-measurable random
variable $\tau_t$ such that $\tau=\tau_t$ on $\{\tau < t\}$.\\
(ii) $\wt Z_\tau=1$ on $\{\tau <\infty\}$.\\
(iii) There exists an optional set $\Lambda$ such that $\tau(\omega)=\sup\{t: (\omega, t)\in \Lambda\}$ on $\{\tau<\infty\}$.\\
(iv) $A^{o}_t=A^{o }_{t\wedge \tau}$. \eethe
\proof Equivalence
between conditions $(i), (ii)$ and $(iii)$ is stated in Theorem
(5,1) from \cite{jeu}. Implication $(i)\Rightarrow (iv)$ comes
from analogous arguments as in \cite{azema}. To finish the proof, we show implication
$(iv)\Rightarrow(iii)$. Let $\Lambda $ be the support of the
measure $dA^o $, i.e.,
\[
\Lambda=\{(\omega, t) | \ \forall \varepsilon>0 \ A^o_t(\omega)>A^o_{t-\varepsilon}(\omega) \}.
\]
The set $\Lambda$ is optional since $A^o$ is an optional process.
Then, $\graph{\tau} \subset \Lambda$ and $A^{o}_t=A^{o }_{t\wedge
\tau}$ imply that indeed $\tau$ is the end of $\Lambda$ on
$\{\tau<\infty\}$. \finproof

\noindent In the case of honest time, any $\ff$-martingale $X$ is
a $\gg$-semimartingale with {(predictable)} decomposition  \cite[Prop. (5,10)]{jeu}
  \be \label{after}X_t= \wh X_t +\int_0^{t\wedge \tau}  \frac {d\cro
{X,m}^\ff_s}{Z_{s-} }-\int _{t\wedge\tau}^t \frac
{d\cro{X,m}^\ff_s}{1-Z_{s-} } \,, \ee where $ \wh  X$
 is a $\gg$-local martingale.

\noindent
 We would like to emphasize the role of $\tilde Z$. As we shall see, this process will be important to prove the existence of arbitrage opportunities.  We give  also a  simple characterisation of honest times avoiding $\ff$- stopping times.
  \bl
\label{stricthonest} A random time $\tau$ is an
honest time and avoids $\ff$-stopping times if and only if $Z
_\tau=1$ a.s. on $(\tau<\infty)$. \el
\proof Assume that $\tau$ is
an  honest time avoiding $\ff$-stopping times. The honesty, by Theorem \ref{honest}, implies
that $\tilde Z_\tau=1$ and the avoiding
property implies the continuity of $A^o$ since for each $\ff$-stopping time $T$, $\E(\Delta A^o_T)=\P(\tau=T<\infty)=0$. Then, the relation $\tilde Z= Z+ \Delta A^{o }$ leads to the result.\\
Assume now  that $Z _\tau=1$ on the set $\{\tau<\infty\}$. Then, on $\{\tau<\infty\}$ we have
$1=Z _\tau\leq \widetilde Z_\tau \leq 1$, so $\widetilde Z_\tau=1$
and $\tau$ is an honest time. Furthermore, as $\Delta
A^o_\tau=\widetilde Z_\tau-Z_\tau=0$, for each $\ff$
stopping time $T$ we have
\begin{align*}
\P(\tau=T<\infty)=\E(\id_{\{\tau=T\}}\id_{\{\Delta A^o_\tau =0 \}}\id_{(T<\infty)})
=\E(\int_0^\infty \id_{\{u=T\}}\id_{\{\Delta A^o_u =0 \}}dA^o_u)=0.
\end{align*}
So $\tau$ avoids $\ff$ stopping times. \finproof

\section{Some particular cases}

\subsection{Immersion assumption, density hypothesis }\label{posdensity}

We recall that the filtration $\ff$ is immersed in $\gg$ under
$\Q$ if any $(\ff,\Q)$-local martingale is a $(\gg,\Q)$-local
martingale.
 \bl If the immersion
property is satisfied under a probability $\Q$ on $\gg$, such that
$S$ is a $(\ff,\Q)$-martingale, all the three concepts of NFLVR, NA and NUPBR hold. \el
\proof Let $S$ be a $(\ff,\Q)$-local
martingale, then it is a $(\gg,\Q)$-local martingale as
well.\finproof

\noindent One says that the random time $\tau$ satisfies the positive density hypothesis if there exists a positive  $\F_t\otimes \mathcal B(\mathbb R^+)$-measurable
function $(\omega, u) \rightarrow \alpha_t( \omega, u)$ which
satisfies:  for any Borel bounded function $\varphi$,  $$ \E(
\varphi(\tau)\vert \F_t)=\int _{\R_+} \varphi (u)\alpha_t(u)
f(u)du
 ,  \,\quad \P-a.s.$$ where $f$ is the density function of $\tau$. In other terms, the conditional distribution of $\tau$
is   characterized by the survival probability  defined by
$$G_t(\theta):= \mathbb P (\tau >\theta \vert \F_t)= \int
_\theta^\infty \alpha_t(u)f(u)du \,.$$ In that case Hypothesis
($\cH^\prime$) is satisfied (see  \cite{am:th} or  \cite{gp:it}).}
\bl If $S$ is a $(\P,\ff)$-martingale and if the conditional law
of $\tau$ with respect to $\ff$ satisfies the positive density
hypothesis then NFLVR holds for $\gg$. Thus both NA and NUPBR hold for $\gg$ as well. \el
\proof
Indeed, under the positive density hypothesis, it can be proved
(see
   Amendinger's thesis
 \cite{am:th} and  Grorud and Pontier \cite{gp:it}), that the
 probability $\P^*$, defined on $\ff \vee \sigma (\tau)$ as
  $$d\P^*\vert_{\F_t \vee \sigma (\tau)}=
\frac{1}{\alpha_t(\tau)}d\P\vert_{\F_t \vee \sigma (\tau)}$$ satisfies the following assertions

(i) Under $\P^*$, $\tau$ is independent from $\F_t$ for any $t$

(ii) $\P^* \vert_{\F_t}=\P \vert_{\F_t}$

(iii) $\P^* \vert_{\sigma(\tau)}=\P \vert_{\sigma(\tau)}$

\noindent Note that immersion is satisfied under ${\mathbb P}^*$. It is
now obvious that, if $S$ is a $({\mathbb P},{\mathbb F})$-martingale, NFLVR holds in
the enlarged
filtration ${\mathbb F} \vee \sigma (\tau)$, hence in $ {\mathbb G}$. Indeed, the $({\mathbb F},{\mathbb P})$-martingale $S$ is
- using the
independence property - an $({\mathbb F}^\tau,{\mathbb P}^*)$-martingale, so that $S$, being $ {\mathbb G}$ adapted, is a $( {\mathbb G},{\mathbb P}^*)$-martingale and ${\mathbb P}^*$ is an equivalent martingale measure. If $S$ is only a local martingale, then one proceeds as follows\footnote{This proof was given to us by C. Fontana}
Let $\{\tau_n\}_{n\in {\mathbb N}}$ be an ${\mathbb F}$-localizing sequence for $S$, meaning that $S^{\tau_n} $ is a $({\mathbb P},{\mathbb F})$ martingale, for every $n\in{\mathbb N}$. Since ${\mathbb P}^*|_{\F_{\infty}}={\mathbb P}|_{\F_{\infty}}$ and  ${\mathbb F}$ is immersed in $ {\mathbb G}$ under ${\mathbb P}^*$, it holds that $S^{\tau_n}$ is a $({\mathbb P}^*, {\mathbb G})$ martingale. Moreover, since $S^{\tau_n}$ is ${\mathbb F}$-adapted, we also have $S^{\tau_n}$ is a $({\mathbb P}^*,{\mathbb F})$ martingale. Finally, the sequence $\{\tau_n\}_{n\in{\mathbb N}}$ is localizing w.r.t. both $({\mathbb P}^*,{\mathbb F})$ and $({\mathbb P}^*, {\mathbb G})$, thus implying that $S$ is a $({\mathbb P}^*,{\mathbb F})$ local martingale. \finproof

\section{Classical arbitrages for a class of honest times}\label{saht}
Herein, we generalize the results obtained in \cite{fjs} -- which
are established for honest times  avoiding  $\ff$-stopping  times
in a complete market with continuous  filtration -- to any
complete market and to a much more broader class of honest times
that will be defined below. Throughout this section, we denote by
${\cal T}_{s}$ the set of all $\mathbb F$-stopping times, ${\cal
T}_h$ the set of all $\mathbb F$-honest times, and ${\cal R}$ the
set of random times given by
 \begin{equation}\label{Rset}
 {\cal R}:=\Bigl\{ \tau\ \ \mbox{random time}\ \big|\ \exists \ \Gamma\in{\cal A}\textrm{ and } T \in {\cal T}_s \textrm{ such that} \ \tau=T\ind_{\Gamma}+\infty \ind _{{\Gamma}^c}\Bigr\},\end{equation}
 { \bp  The following inclusions hold
 \begin{equation}\label{inclusions}
 {\cal T}_{s}\subset {\cal R}\subset {\cal T}_h.
 \end{equation}
 \ep}
\noindent \proof The first inclusion is clear. For  the inclusion ${\cal R}\subset {\cal T}_h$, we give, for ease of the reader two different proofs. Let us take $\tau \in \mathcal R$.\\
1) On $(\tau<t)=(T<t)\cap {\Gamma}$, we have $\tau=T\land t$ and $T\land t$ is $\F_t$-measurable. Thus, $\tau$ is an honest time.\\
2)
We want to show that on $(\tau<\infty)$, $\widetilde Z_\tau=1$. Indeed, $\widetilde Z_t=\id_{(T\geq t)}\P(\Gamma|\F_t)+\P(\Gamma^c|\F_t)$, so that
\[
\id_{(\tau<\infty)}\widetilde Z_\tau=\id_{\Gamma}\id_{(T<\infty)}\widetilde Z_T
=\id_\Gamma\id_{(T<\infty)}(\id_{(T\geq T)}\P(\Gamma|\F_T)+\P(\Gamma^c|\F_T))=\id_{\Gamma}\id_{(T<\infty)}=\id_{(\tau<\infty)}.
\]
\noindent This proves that $\tau$ is an honest time.\finproof

\noindent The following theorem  represents our principal result
in the general framework.

\bethe\label{sa} Assume that $(S ,\ff)$  is a  complete market and let $\varphi$ be an $\ff$-predictable process satisfying $m=1+\varphi\centerdot S$.
Then the following assertions hold.
\begin{enumerate}
\item[(a)] If $\tau$ is an honest time, and $\tau \not\in {\cal R}$, then $\gg$-predictable process $\varphi^b=\varphi\id_{\Rbrack0,\tau\Lbrack}$ is a classical arbitrage strategy in the market ''before $\tau$'', i.e., in $(S^\tau,\gg)$.
\item[(b)] If $\tau$ is an honest time, which is not an $\ff$-stopping time, and if $\{\tau=\infty\}\in \F_\infty$, then the $\gg$-predictable process $\varphi^a=-\varphi \id_{\Lbrack\tau, \nu\Lbrack}$, with $\gg$-stopping time defined as
\begin{equation}\label{nu}
\nu:=\inf\{t>\tau: \wt Z_t\leq \frac{1-\Delta A^o_\tau}{2}\},\end{equation}
 is a classical arbitrage strategy in the market  ''after $\tau$'', i.e., in $(S-S^\tau,\gg)$.
\end{enumerate}
 \eethe
\proof  {(a)} From $m= \wt Z+ A^o_-$ and $\wt Z_\tau=1$, we deduce that $m_\tau\geq 1$. Since $\tau\notin \mathcal R$, one has $\P(m_\tau>1)=\P(  A^o_{\tau-}>0)>0$.  Then, by Lemma  \ref{lemmem}, process $\varphi^b=\varphi\id_{\Rbrack 0, \tau \Lbrack}$ is an arbitrage strategy in $(S^\tau, \gg)$.\\
(b) From $m= Z+A^o$  and Theorem \ref{honest} (iv), one obtains that, for $t>\tau$,  $ m_t-m_\tau =Z_t-Z_\tau  \geq -1$.
On the other hand, using $m= \wt Z+A^o_-$, one obtains that, for $t>\tau$,  $ m_t-m_\tau =\wt Z_t-1 +\Delta A^o_\tau$.
Assumption $\{\tau=\infty\}\in \F_\infty$ ensures that  $\wt Z_\infty= \id_{\{\tau=\infty\}}$ and in particular $\{\tau<\infty\}\subset\{ \wt Z_{\infty }=0\}$.
So, $\gg$-stopping time $\nu$ defined in \eqref{nu} satisfies $\{\nu<\infty\}=\{\tau <\infty\}$.
Then, $$m_\nu-m_\tau=\wt Z_\nu-1 +\Delta A^o_\tau\leq \frac{\Delta A^o_\tau-1}{2}\leq 0,$$
and, as $\tau$ is not  an $\ff$-stopping time,
$$\P(m_\nu-m_\tau<0)=\P(\Delta A^o_\tau <1)>0.$$
Hence $-\int_\tau ^{t\land\nu} \varphi_s dS_s =m_{\tau\land t}-m_{t\land \nu}$  is the value of an admissible  self-financing strategy $\varphi^a=-\varphi\id_{\Lbrack \tau, \nu\Lbrack}$ with initial value 0 and terminal value $m_\tau-m_\nu\geq 0$ satisfying $\P(m_\tau-m_\nu>0)>0$. This ends the proof of the theorem.
 \finproof

\brem We recall that if $\tau$ is a finite  honest  time (hence
$\F_\infty$-measurable) and is not an $\ff$-stopping time, then
the density hypothesis is not satisfied and immersion does not
hold. Indeed:
\begin{enumerate}
\item[(i)] Density hypothesis would  hold if, under some equivalent
probability measure, $\tau$ would be independent from
$\F_\infty$.
\item[(ii)]The immersion  property is equivalent to
$\P(\tau >t \vert \F_t)=\P(\tau >t \vert \F_\infty)$ which, for a
finite honest time is $\ind_{\tau>t}$. Then, one should have
$\P(\tau >t \vert \F_t)=\ind_{\tau>t}$ and $\tau$ would be a
stopping time.
\end{enumerate}
\erem
\brem The completeness of the market is
an obvious condition to conclude. See \cite{fjs} for a counter
example.\erem

\noindent In the following two subsections we explore several examples of honest times.
Each of them is defined as an end of optional set, so by Theorem \ref{honest} (iii), is indeed an honest time.

\subsection{Classical arbitrage opportunities in a Brownian filtration}

In this subsection, we develop practical market models $S$ and
honest times $\tau$ within the Brownian filtration for which
one can  compute explicitly the  arbitrage opportunities for both
before and after $\tau$. For other examples of honest times,
and associated classical arbitrages  we refer the reader to
\cite{fjs} (note that the arbitrages constructed in that paper are
different from our arbitrages). Throughout this subsection, we
assume given a one-dimensional Brownian motion $W$ and $\ff$ is
its augmented natural filtration. The market model is represented
by the bank account whose process is the constant one and one
stock whose price process is given by
$$S_t=  \exp (\sigma
W_t-\frac 12 \sigma^2t),\ \ \ \ \ \sigma>0\  \mbox{given} .$$
It is worth mentioning that in this context of Brownian filtration, for any process $V$ with locally integrable variation, its $\ff$-dual optional projection is equal to its $\ff$-dual predictable projection, i.e., $V^{o,\ff}=V^{p,\ff}$.

\subsubsection{Last passage time  at a given level}

\bp Consider the following random times
 $$\tau:=\sup\{t\,:\,
S_t=a\}\ \ \ \ \ \mbox{and}\ \ \ \ \nu:=\inf\{t>\tau\ \big|\ \ S_t\leq {{a}\over{2}}\},$$ where $0<a< 1$. Then, the following assertions hold.\\
(a) The model "before $\tau$" $(S^{\tau},\gg)$
 admits a classical arbitrage opportunity  given by the $\gg$-predictable process $$ \varphi^b=\frac{1} {a}\ind_{\{S<a\}} I_{\Lbrack0,\tau\Lbrack} .$$
(b) The model "after $\tau$" $(S-S^{\tau},\gg)$ admits a classical arbitrage
opportunity given by $\gg$-predictable process
$$\varphi^a= - \frac{1} {a}\ind_{\{S<a\}} I_{\Lbrack\tau,\nu\Lbrack} .$$
 \ep
\proof  Since $\tau\in {\cal T}_h \backslash \cal R$ we make a use of Theorem \ref{sa}. We compute the predictable process
$\varphi$ such that $m=1+\varphi \centerdot S$.   To this end, we
calculate $Z$ as follows. Using \cite[exercise 1.2.3.10]{3M}, we
derive
$$1-Z_t:=\P\left(\tau \leq t\vert \F_t\right)= \P\left(\sup_{t <u }S_u \leq a\vert \F_t\right)=\P\left(\sup_{u } \widetilde S_u \leq \frac a {S_t} \vert \F_t\right) = \Phi\left( \frac a {S_t}\right) $$
 where $\widetilde  S_u=\exp ( \sigma \widetilde W_u-\frac 12 \sigma ^2u)$, $\widetilde W$ independent of $\F_t$ and
 $\Phi(x)= \P\left(\sup_{u } \widetilde S_u \leq x\right)= \P(\frac 1 U \leq x) = \P(\frac 1 x \leq U)= (1-\frac 1x)^+$,
where $U$ is a random variable with uniform law. Thus we get $Z_t= 1- (1-\frac {S_t}a)^+$ (in particular $Z_\tau=\wt Z_\tau=1$), and  $$dZ_t =
 \ind_{\{S_t<a\}} \frac 1 a dS_t  -\frac 1{2a} d\ell^a_t$$
where $\ell^a$ is the local time of the $S$ at the level $a$ (see page 252 of He et al.  \cite{yanbook} for the definition of the local time).
Therefore,  we deduce that
$$m=1+\varphi\is S.$$
Note that $\nu:=\inf\{t>\tau\ \big|\ \ S_t\leq {{a}\over{2}}\}=\inf\{t>\tau \ | \ 1-(1-\frac{S_t}{a})^+\leq \frac{1}{2}\}$, so $\nu$ coincides with \eqref{nu}.
Theorem \ref{sa} ends the proof of the proposition. \finproof

\subsubsection{Last passage time at a level before maturity }

 \noindent Our second example of random time, in this subsection, takes into account finite horizon.
 In this example, we introduce the following notation
 \be\label{HVprocess}
 H(z,y,s):=e^{-z y} \N \left ( \frac{z s -y}{\sqrt
s}\right) +e^{z y}\N \left ( \frac{-z s -y}{\sqrt
s}\right),\ee
where $\N(x)$ is the cumulative distribution function of the standard normal distribution.

\bp
Consider the following random time (an honest time)
$$\tau_1:=\sup\{t \leq 1\,:\, S_t= b\}   $$
where $b$ is a positive real number, $0<b<1$ .
Let $V$ and $\beta$ be given by
\[
V_t:=\alpha -\gamma t - W_t\,\ \ \mbox{with}\ \alpha=\frac{\ln b}{\sigma}\ \mbox{and} \ \gamma=-\frac{\sigma}{2}
\]
$$
\beta_t:=e^{ \gamma V_t}\left(\gamma H(\gamma, \vert V_t\vert
, 1-t) - \mbox{sgn} (V_t) \,H^\prime _x( \gamma, \vert  V_t\vert ,
1-t)\right),$$ with $H$ defined in (\ref{HVprocess}), and let
$\nu$ be as in (\ref{nu}).
Then, the following assertions hold.\\
(a) The model "before $\tau_1$" $(S^{\tau_1},\gg)$
 admits a classical arbitrage opportunity  given by the $\gg$-predictable process $$ \varphi^b:= \frac{1}{\sigma S_t} \beta_t I_{\Rbrack 0,\tau_1\Lbrack}.$$
(b) The model "after $\tau_1$" $(S-S^{\tau_1},\gg)$ admits a classical arbitrage
opportunity given by $\gg$-predictable process
$$\varphi^a := -\frac{1}{\sigma S_t} \beta_t I_{\Lbrack\tau_1,\nu\Lbrack} .$$
 \ep
\proof The proof of this proposition follows from Theorem \ref{sa} as long as we can write the martingale $m$ as an integral stochastic with respect to $S$. This is the main focus of the remaining part of this proof.
By Theorem \ref{honest} (iii), the time $\tau_1$ is honest and finite.
Honest time $\tau_1$ can be seen as \beq
\tau_1  &= &\sup \, \{ t\leq  1 \,:\,\gamma t + W_t = \alpha \}
=\sup \, \{ t\leq  1 \,:\,V_t = 0 \}.
 \eeq
Setting  $T_0(V)=\inf\{t\,:\, V_t=0\}$, we
obtain, using standard computations (see \cite{3M} p. 145-148)  $$1-Z_t=\P(
\tau_1 \leq t\vert \F_t)= (1 - e^{ \gamma V_t }H(\gamma, \vert V_t
\vert , 1-t)) \ind_{\{ T_0 (V) \leq t \leq 1\}}+\id_{\{t>1\}},\,
$$   where $H$ is given in (\ref{HVprocess}). In particular
$Z_\tau=\wt Z_\tau =1$.
 Using It\^o's lemma, we obtain  the decomposition of $1 -
e^{ \gamma V_t }H(\gamma, \vert V_t \vert , 1-t)$ as a
semimartingale. The martingale part of $Z$ is  given
  by  $dm_t=\beta_tdW_t= \frac{1}{\sigma S_t} \beta_tdS_t$, which ends the proof.
\finproof

\subsection{Arbitrage opportunities in a Poisson filtration }\label{poisson}

Throughout this subsection, we suppose given a Poisson process  $N$, with intensity rate $\lambda>0$, and natural filtration  $\ff$. The stock price process is given by
\begin{equation}\label{modelpoisson2}
dS_t=S_{t-}\psi  dM_t,\, \ \ S_0=1,\ \ \ \ \ M_t:=N_t-\lambda t,\end{equation}
 or equivalently $S_t= \exp (-\lambda \psi t +\ln (1+\psi) N_t)$, where $\psi>-1$. In what follows, we introduce the notation
\begin{equation}\label{modelpoisson2}
 \alpha:=\ln (1+\psi),\ \  \ \ \mu:= \frac{\lambda \psi}{\ln (1+\psi)}\ \ \mbox{and}\ \  Y_t:=  \mu t -  N_t,\end{equation}
 so that $S_t=\exp (-\ln (1+\psi) Y_t)$. We associate to the process $Y$ its ruin probability, denoted by
$\Psi(x)$   given by  \be \label{Psi0}\Psi(x)=\P(T^x<\infty)
,\quad \mbox{ with} \quad
 T^x=\inf\{t: x+Y_t <0\}\,\ \ \ \ \mbox{and}\ \ \ \ x\geq 0.\ee
\noindent Below, we describe our first example of honest time and the associated arbitrage opportunity.

\subsubsection{Last passage time at a given level}
\bp Suppose that $\psi>0$ and let $\varphi$ be
  $$\varphi:={{\Psi(Y_{-}-a-1)\ind_{\{Y_{-} \geq a+1\}}-\Psi(Y_{-}-a)\ind_{\{Y_{-} \geq a\}}+\ind_{\{Y_{-} < a+1\}} -\ind_{\{Y_{-} < a\}}}\over{\psi S_{-}}}.$$
For $0<b<1$, consider the following random time
\begin{equation}\label{suppoi}
\tau:=  \sup \{t:\,S_t \geq b\}=\sup \{t: Y_t \leq a \},
\end{equation}
with $a:= -\frac{1}{\alpha}\ln b$.
Then the following assertions hold.\\
(a) The model "before $\tau$" $(S^{\tau},\gg)$
 admits a classical arbitrage opportunity  given by the $\gg$-predictable process $ \varphi^b:= \varphi I_{\Rbrack 0,\tau_1\Lbrack}.$\\
(b) The model "after $\tau$" $(S-S^{\tau},\gg)$ admits a classical arbitrage
opportunity given by $\gg$-predictable process
$\varphi^a := -\varphi I_{\Lbrack\tau_1,\nu\Lbrack},$ with $\nu$ as in \eqref{nu}.
\ep
\proof
   Since $\psi>0$, one has   $ \mu >\lambda$   so  that $Y$ goes to $+\infty$ as $t$ goes to infinity, and $\tau$ is finite.
  The supermartingale $Z$
associated with the   time $\tau$ is
$$Z_t=
\P(\tau>t|\F_t)
 =\Psi(Y_t-a)\ind_{\{Y_t \geq a\}} + \ind_{\{Y_t < a\}}= 1+\ind_{\{Y_t \geq a\}}\left(\Psi(Y_t-a)-1\right),
$$
where $\Psi$ is defined in (\ref{Psi0}) (see \cite{thin} for more details on this example).
We set $\theta = \displaystyle\frac{\mu}{\lambda   }-1$,  and deduce that $\Psi (0)= (1+\theta)^{-1}$ (see  \cite{as}).   Define $\vartheta_1=\inf\{t> 0: Y_t = a\}$
and then, for each $n>1$, $\vartheta_n=\inf\{t> \vartheta_{n-1}:
Y_t = a\}$. It can be proved that the times $\vartheta_n$ are
predictable $\ff$-stopping times (\cite{thin}).
The $\ff$-dual optional projection  $A^{o}$  of the process $\ind_{\Rbrack\tau, \infty\Lbrack}$ equals
$$A^{o}= \frac{\theta}{1+\theta} \sum_n\ind_{\Rbrack\vartheta_n,\infty\Lbrack}.$$
 Indeed, for any $\ff$-optional process $U$ we have
$$ \E( U_\tau)  = \E(\sum \ind_{\{\tau = \vartheta_n\}}
U_{\vartheta_n})=\E(\sum \E(\ind_{\{\tau = \vartheta_n\}} \vert
\F_{\vartheta_n}) U_{\vartheta_n})$$ and $\E(\ind_{\{\tau =
\vartheta_n\}} \vert \F_{\vartheta_n})= \P( T^0=\infty)= 1-\Psi(0) =
1-\frac{1}{1+\theta}$.

\noindent As a result the process $A^o$ is predictable, and hence
$Z=m-A^o$ is the Doob-Meyer decomposition of $Z$. Thus we can get
$$
\Delta m= Z-\ ^{p}Z$$  where  $^{p}Z $ is the $\ff$-predictable
projection of $Z$ \footnote{Note that here we talk about predictable projection and not about dual predictable projection.}. To calculate $^{p}Z$, we   write the process
$Z$ in a more adequate form. To this end, we first remark that
$$\ind_{\{Y \geq a\}}=\ind_{\{Y_{-} \geq a+1\}} \Delta N+(1-\Delta N)\ind_{\{Y_{-} \geq a\}}\ \ \ \ \mbox{and}\ \ \  \ind_{\{Y < a\}}=\ind_{\{Y_{-} < a+1\}} \Delta N+(1-\Delta N)\ind_{\{Y_{-} <a\}}.$$ Then, we obtain
$$\begin{array}{lll}
\Delta m=\left(\Psi(Y_{-}-a-1)\ind_{\{Y_{-} \geq a+1\}}-\Psi(Y_{-}-a)\ind_{\{Y_{-} \geq a\}}+\ind_{\{Y_{-} < a+1\}} -\ind_{\{Y_{-} < a\}}\right)\Delta N \\
\\
\hskip 1cm =\psi S_{-}\varphi\Delta M=\varphi\Delta
S.\end{array}$$
Since the two martingales $m$ and $S$ are
discontinuous, we deduce that $m=1+\varphi\centerdot S$.
Therefore, the proposition follows from Theorem \ref{sa}.
 \finproof

\subsubsection{Time of supremum on fixed time horizon }
 \noindent The second example requires the following notations
\be \label{Psi1/Phihat} S^*_t:= \sup _{s\leq t}S_ s,\ \ \ \Psi(x,
t):= \P( S_t^* > x) ,\ \ \ \ \wh \Phi( t):=\P(\sup_{s<t}S_s\leq
1),\ \ \ \wt \Phi(x, t):=\P(\sup_{s<t}S_s<x)\ee

\bp
 Consider the random time $\tau$ defined by \begin{equation}\label{honnete1}
\tau = \sup \{t\leq 1: S_t=  S_t^* \},
\end{equation}  {where $S^*_t=\sup_{s\leq t}S_s$}.
Then,  the following assertions hold.\\
a) For $\psi >0$, define the $\gg$-predictable process $\varphi$ as
$$\begin{array}{llll}
\varphi_t:=\ind_{\{ t<1\}}\left[\Psi\left(\max({{S^*_{t-}}\over{S_{t-}(1+\psi)}},1),1-t\right)-\Psi\left({{S^*_{t-}}\over{S_{t-}}},1-t\right)\right]+\id_{\{S^*_{t-}<S_{t-}(1+\psi)\}}\;\wh \Phi(1-t)\\
\\
\hskip 1cm +\left[\id_{\{\max(S_{1-}^*, S_{1-}(1+\psi))=S_0\}}-\id_{\{\max(S_{1-}^*, S_{1-})=S_0\}}\right]\ind_{\{ t=1\}}\,\end{array}.$$
Then, $\varphi^b:=\varphi\id_{\Rbrack0, \tau\Lbrack}$
is an arbitrage opportunity for the model $(S^{\tau},\gg)$,  and $\varphi^a:=-\varphi I_{\Lbrack\tau,\nu\Lbrack}$  is an arbitrage opportunity for the model $(S-S^{\tau},\gg)$.  Here
$\Psi$ and $\wh \Phi$ are defined in
(\ref{Psi1/Phihat}), and $\nu$ is defined similarly as in (\ref{nu}).\\
b)  For $-1<\psi <0$, define the $\gg$-predictable process
$$\varphi_t:=\frac{\psi I_{\{ S^*_t=S_{t-}\}}\wh\Phi({1\over{1+\psi}},1-t)+\Psi({{S^*_t}\over{S_{t-}(1+\psi)}},1-t)
-\Psi({{S^*_t}\over{S_{t-}}},1-t)}{\psi S_{t-}}.$$
Then, $\varphi^b:=\varphi\id_{\Rbrack0, \tau\Lbrack}$
is an arbitrage opportunity for the model $(S^{\tau},\gg)$,  and $\varphi^a:=-\varphi I_{\Lbrack\tau,\nu\Lbrack}$  is an arbitrage opportunity
for the model $(S-S^{\tau},\gg)$.
\ep
\proof Note that, if
$-1<\psi <0$ the process $S^*$ is continuous, $S_\tau< S_\tau^*=\sup_{t\in[0,1]}S_t$ on the set
$(\tau<1)$ and
$S_{\tau-}= S_{\tau-}^*=\sup_{t\in[0,1]}S_t$.
If  $\psi>0$, $S_{\tau-}<S_{\tau-}^*<\sup_{t\in[0,1]}S_t$  on the set $(\tau<1)$.

\noindent Define the sets $(E_n)_{n=0}^\infty$ such that
$E_0=\{\tau=1\}$ and $E_n=\{\tau = T_n\}$ with $n\geq 1$.
The sequence $(E_n)_{n=0}^\infty$ forms a partition of $\Omega$.
Then, $\tau=\id_{E_0}+\sum_{n=1}^\infty T_n\id_{E_n}$.
Note that $\tau$ is not an $\ff$ stopping time since  $E_n \notin \F_{T_n}$ for any $n\geq 1$.

\noindent The supermartingale $Z$ associated with the honest time $\tau$
is
 \begin{align*}
Z_t=\P(\tau>t|\F_t)&= \P(\sup_{s\in(t,1]}S_s>
\sup_{s\in[0,t]}S_s|\F_t) =\P(\sup_{s\in[0,1-t]}\widehat S_s>
\frac{  S_t^*}{S_t}|\F_t) =\id_{(t<1)}\Psi(\frac{  S_t^*}{S_t},
1-t),
\end{align*}
with $\widehat S$ an independent copy of $S$ and $\Psi(x,t)$ is given by (\ref{Psi1/Phihat}).

\noindent As $\{\tau=T_n\}\subset\{\tau\leq T_n\}\subset \{Z_{T_n}<1\}$,
we have
\begin{align*}
Z_\tau=\id_{\{\tau=1\}}Z_{1}+\sum_{n=1}^\infty
\id_{\{\tau=T_n\}}Z_{T_n}<1,\ \ \ \mbox{and}\ \ \{\widetilde Z=0<Z_{-}\}=\emptyset.
\end{align*}
\noindent In the following we will prove assertion a). Thus, we suppose that $\psi>0$, and we calculate
\begin{align*}
A^o_t&=\P(\tau=1|\F_1)\id_{\{t\geq 1\}}+\sum_n \P(\tau=T_n|\F_{T_n})\id_{\{t\geq T_n\}}\\
&=\id_{\{S_1^*=S_0\}}\id_{\{t\geq  1\}}+\sum_n \id_{\{T_n<1\}}\id_{\{S^*_{T_n-}<S_{T_n}\}}\P(\sup_{s\in [T_n,1[}S_s\leq S_{T_n}|\F_{T_n})\id_{\{t\geq T_n\}}\\
&=\id_{\{S_1^*=S_0\}}\id_{\{t\geq 1\}}+\sum_n \id_{\{T_n<1\}}\id_{\{S^*_{T_n-}<S_{T_n-}(1+ \psi)\}}\;\wh\Phi(1-T_n)\id_{\{t\geq T_n\} },
\end{align*}
with $\wh\Phi$ is given by (\ref{Psi1/Phihat}). As
 before, we write
\begin{align*}
A^o_t&=\id_{\{S_1^*=S_0\}}\id_{\{t\geq 1\}}+\sum_{s\leq t} \id_{\{s<1\}}\id_{\{S^*_{s-}<S_{s-}(1+\psi)\}}\;\wh \Phi(1-s)\Delta  N_s\\
&=\id_{\{S_1^*=S_0\}}\id_{\{t\geq 1\}}+\int_0^{t\land 1}\id_{\{S^*_{s-}<S_{s-}(1+\psi)\}}\;\wh \Phi(1-s)\, d M_s
+\lambda \int_0^{t\land 1} \id_{\{S^*_{s-}<S_{s-}(1+\psi)\}}\;\hat \Phi(1-s) d s.
\end{align*}

\noindent Remark that we have
$$
\ind_{\{S_1^*=S_0\}}=\left[\ind_{\{\max(S_{1-}^*, S_{1-}(1+\psi))=S_0\}}-\id_{\{\max(S_{1-}^*, S_{1-})=S_0\}}\right]\Delta M_1+\ind_{\{\max(S_{1-}^*, S_{1-})=S_0\}}.$$
and
$$
\Delta m=\Delta Z+\Delta A^o= Z-\ ^{p}(Z)+\Delta A^o-\ ^{p}(\Delta A^o).$$
Then we re-write the process $Z$ as follows
$$
Z=\ind_{\Rbrack 0, 1\Rbrack}\Psi\left(\max({{S^*_{-}}\over{S_{-}(1+\psi)}},1),1-t\right)\Delta M+(1-\Delta M)I_{\Rbrack 0, 1\Rbrack}\Psi\left({{S^*_{-}}\over{S_{-}}},1-t\right).$$ This implies that $$
Z-\ ^{p}(Z)=\ind_{\Rbrack 0, 1\Rbrack}\left[\Psi\left(\max({{S^*_{-}}\over{S_{-}(1+\psi)}},1),1-t\right)-\Psi\left({{S^*_{-}}\over{S_{-}}},1-t\right)\right]\Delta M.$$
Thus by combining all these remarks, we deduce that
$$\begin{array}{lll}
\Delta m=Z-\ ^{p}(Z)+\Delta A^o-\ ^{p}(\Delta A^o)=\varphi \Delta S.\end{array}
$$
Then, the assertion a) follows immediately from Theorem \ref{sa}.\\
Next, we will prove assertion b). Suppose that $-1<\psi<0$, and we calculate
\begin{align*}
A^o_t&=\P(\tau=1|\F_1)\id_{\{t\geq 1\}}+\sum_n \P(\tau=T_n|\F_{T_n})\id_{\{t\geq T_n\}}\\
&=\id_{\{S_1^*=S_1\}}\id_{\{t\geq 1\}}+\sum_n \id_{\{T_n<1\}}\id_{\{S^*_{T_n}=S_{T_n-}\}}\P(\sup_{s\in [T_n,1[}S_s<S_{T_n-}|\F_{T_n})\id_{\{t\geq T_n\}}\\
&=\id_{\{S_1^*=S_1\}}\id_{t\geq 1\}}+\sum_n \id_{\{T_n<1\}}\id_{\{S^*_{T_n}=S_{T_n-}\}}\wt \Phi(\frac{S_{T_n-}}{S_{T_n}},1-T_n)\id_{\{t\geq T_n\}},
\end{align*}
with $\wt \Phi(x, t)$ is given by (\ref{Psi1/Phihat}). In order
to find  the compensator of $A^o$,   we write
\begin{align*}
A^o_t&=\id_{\{S_1^*=S_1\}}\id_{\{t\geq 1\}}+\sum_{s\leq t} \id_{\{s<1\}}\id_{\{S^*_{s}=S_{s-}\}}\wt \Phi(\frac{1}{  1+\psi },1-s)\,\Delta  N_s\\
&=\id_{\{S_1^*=S_1\}}\id_{\{t\geq 1\}}+\int_0^{t\land 1} \id_{\{S^*_{s}=S_{s-}\}}\wt \Phi(\frac{1}{  1+\psi },1-s)\, d M_s
+\lambda \int_0^{t\land 1} \id_{\{S^*_{s}=S_{s-}\}}\wt \Phi(\frac{ 1}{  1+\psi },1-s)\, d s.
\end{align*}
As a result, due to the continuity of the process $S^*$,  we get
\beq A^o_t-\ ^p(A^o)_t&=&I_{\{
S^*_t=S_{t-}\}}\wt\Phi({1\over{1+\psi}},1-t)\Delta M_t,\\  Z_t-\
^pZ_t&=&\left[\Psi({{S^*_t}\over{S_{t-}(1+\psi)}},1-t)-\Psi({{S^*_t}\over{S_{t-}}},1-t)\right]\Delta
N_t.\eeq This implies that \beq  \Delta m_t&=&Z_t-\ ^pZ_t+A^o_t-\
^p(A^o)_t\\&=&\left\{\psi
I_{\{S^*_t=S_{t-}\}}\wt\Phi({1\over{1+\psi}},1-t)+\Psi\left({{S^*_t}\over{S_{t-}(1+\psi)}},1-t\right)-\Psi\left({{S^*_t}\over{S_{t-}}},1-t\right)\right\}\Delta
N_t.\eeq
 Since $m$ and $S$
are pure discontinuous $\ff$-local martingales, we conclude that $m$ can
be written in the form of
$$
m=m_0+\varphi\cdot S,$$
 and the proof of the assertion b) follows immediately from Theorem \ref{sa}. This ends the proof of the proposition.\finproof

\subsubsection{Time of overall supremum}\label{sup}
\noindent Below, we will present our last example of this subsection. The analysis of this example is based on the following three functions.

\be\label{threefunctionPhi} \Psi(x )= \P( S ^* > x)=\P( \sup_sS_s
> x),\ \ \ \  \wh\Phi =\P(\sup_{s}S_s\leq 1),\ \mbox{and}\ \ \ \wt
\Phi(x)=\P(\sup_{s}S_s<x). \ee \bp Consider the random time $\tau$
given by
\begin{equation} \label{honnete2}
\tau = \sup \{t : S_t=  S_t^* \}.
\end{equation}
Then,  the following assertions hold.\\
a) For $\psi >0$, define the $\gg$-predictable process $\varphi$ as
$$
\varphi_t:={{\id_{\{S^*_{t-}<S_{t-}(1+\psi)\}}\wh\Phi+\Psi\left(\max({{S^*_{t-}}\over{S_{t-}(1+\psi)}},1\right)-\Psi({{S^*_{t-}}\over{S_t{-}}})}\over{S_{t-}\psi}}.$$
Then, $\varphi^b:=\varphi\id_{\Rbrack0, \tau\Lbrack}$
is an arbitrage opportunity for the model $(S^{\tau},\gg)$,  and $\varphi^a:=-\varphi I_{\Lbrack\tau,\nu\Lbrack}$  is an arbitrage opportunity for the model $(S-S^{\tau},\gg)$.
Here $\Psi$ and $\wh \Phi$ are defined in
(\ref{threefunctionPhi}), and $\nu$ is defined in similar way as in (\ref{nu}).\\
b) For $-1<\psi<0$, define the $\gg$-predictable process $\varphi$ as
$$
\varphi:={{\Psi(\frac{  S^*}{S_{-}(1+\psi)})-\Psi(\frac{  S^*}{S_{-}})+\id_{\{S^*=S_{-}\}}\wt \Phi(\frac{1}{  1+\psi } )\psi}\over{\psi S_{-}}}.$$
Then, $\varphi^b:=\varphi\id_{\Rbrack0, \tau\Lbrack}$
is an arbitrage opportunity for the model $(S^{\tau},\gg)$,  and $\varphi^a:=-\varphi I_{\Lbrack\tau,\nu\Lbrack}$  is an arbitrage opportunity for the model $(S-S^{\tau},\gg)$.
Here again $\nu$ is defined as in (\ref{nu}).
\ep
\proof Let us note that $\tau$ is  {finite}  and, as before,    if    $-1<\psi <0$,
$S_\tau< S_\tau^*=\sup_{t }S_t$ and $S^*$ is continuous
and if   $\psi >0$, $S_\tau=
S_\tau^*=\sup_{t }S_t$.
\\
The supermartingale $Z$ associated with the honest time $\tau$
is
\begin{align*}
Z_t=\P(\tau>t|\F_t)&= \P(\sup_{s\in(t,\infty]}S_s >
\sup_{s\in[0,t]}S_s|\F_t) =\P(\sup_{s\in[0,\infty]}\widehat
S_s > \frac{   S_t^*}{S_t}|\F_t) = \Psi(\frac{  S_t^*}{S_t}),
\end{align*} with $\widehat S$ an independent copy of $S$ and  $\Psi $ is given by (\ref{threefunctionPhi}). As a result, we deduce that
$ Z_{\tau}<1$.
\noindent In the following, we will prove assertion a). We suppose
that $\psi>0$, denoting by $(T_n)_n$ the sequence of jumps of the
Poisson process $N$,  we derive
\begin{align*}
A^o_t&=\sum_n \P(\tau=T_n|\F_{T_n})\id_{\{t\geq T_n\}}
 =\sum_n \id_{\{S^*_{T_n-}<S_{T_n}\}}\P(\sup_{s\geq T_n}S_s\leq S_{T_n}|\F_{T_n})\id_{\{t\geq T_n\}} \\
&=\sum_n \id_{\{S^*_{T_n-}<S_{T_n-}(1+\psi)\}}\wh \Phi
\id_{\{t\geq T_n\}},
\end{align*}
with $\wh\Phi =\P(\sup_{s}S_s\leq 1)$ given by (\ref{threefunctionPhi}).

\noindent We continue to find compensator of $A^o$
\begin{align*}
A^o_t&=\sum_{s\leq t} \id_{\{S^*_{s-}<S_{s-}(1+\psi)\}}\wh \Phi \Delta  N_s\\
&=\int_0^{t} \id_{\{S^*_{s-}<S_{s-}(1+\psi)\}}\wh\Phi d M_s
+\lambda \int_0^{t}  \id_{\{S^*_{s-}<S_{s-}(1+\psi)\}}\wh \Phi d s.
\end{align*}

\noindent Now as we did for the previous propositions, we calculate the jumps of $m$. To this end, we re-write $Z$ as follows
$$
Z=\left[\Psi\left(\max({{S^*_{-}}\over{S_{-}(1+\psi)}},1)\right)-\Psi({{S^*_{-}}\over{S_{-}}})\right]\Delta M+ \Psi({{S^*_{-}}\over{S_{-}}}).$$ This implies that
$$
Z-\,^{p}Z=\left[\Psi\left(\max({{S^*_{-}}\over{S_{-}(1+\psi)}},1)\right)-\Psi({{S^*_{-}}\over{S_{-}}})\right]\Delta
M.$$ Hence, we derive
$$
\Delta m= \left[\id_{\{S^*_{s-}<S_{s-}(1+\psi)\}}\wh\Phi+\Psi\left(\max({{S^*_{-}}\over{S_{-}(1+\psi)}},1)\right)-\Psi({{S^*_{-}}\over{S_{-}}})\right]\Delta M.$$ \noindent Since both martingales $m$ and $M$ are
purely discontinuous, we deduce that $m=m_0+\varphi\centerdot S$. Then, the proposition follows immediately from Theorem \ref{sa}.\\
\noindent In the following, we will prove assertion b). To this end, we suppose that $\psi<0$, and we calculate
\begin{align*}
A^o_t&=\sum_n \P(\tau=T_n|\F_{T_n})\id_{\{t\geq T_n\}}
 =\sum_n \id_{\{S^*_{T_n}=S_{T_n-}\}}\P(\sup_{s\geq T_n}S_s<S_{T_n-}|\F_{T_n})\id_{\{t\geq T_n\}}\\
&=\sum_n \id_{\{S^*_{T_n}=S_{T_n-}\}}\wt \Phi(\frac{S_{T_n-}}{S_{T_n}})\id_{\{t\geq T_n\}},
\end{align*}
with $\wt \Phi(x)=\P(\sup_{s}S_s<x)$. Therefore,
\begin{align*}
A^o_t&=\sum_{s\leq t} \id_{\{S^*_{s}=S_{s-}\}}\wt \Phi(\frac{1}{ 1+\psi } )\Delta  N_s\\
&=\int_0^{t} \id_{\{S^*_{s}=S_{s-}\}}\wt \Phi(\frac{1}{  1+\psi } ) d M_s
+\lambda \int_0^{t} \id_{\{S^*_{s}=S_{s-}\}}\wt \Phi(\frac{1}{  1+\psi } ) d s.
\end{align*}
Since in the case of $\psi<0$, the process $S^*$ is continuous, we obtain
$$
Z-\ ^pZ=\left[\Psi(\frac{  S^*}{S_{-}(1+\psi)})-\Psi(\frac{
S^*}{S_{-}})\right]\Delta N,\ \ \ \ A^o-\
^{p}(A^o)=\id_{\{S^*=S_{-}\}}\wt \Phi(\frac{1}{  1+\psi } ) \Delta
M.$$ Therefore, we conclude that
$$
\Delta m=Z-\ ^pZ+
A^o-\ ^{p}(A^o)=\left\{\Psi(\frac{
S^*}{S_{-}(1+\psi)})-\Psi(\frac{
S^*}{S_{-}})+\id_{\{S^*=S_{-}\}}\wt \Phi(\frac{1}{  1+\psi }
)\psi\right\}\Delta N.$$ This implies that the martingale $m$ has
the form of $m=1+\varphi\cdot S$, and assertion b) follows
immediately from Theorem \ref{sa}, and the proof of the
proposition is completed.\finproof

\section{Arbitrage opportunities for non-honest random times}
This section is our second main part of the corp of the paper. Herein, we develop a number of practical examples of market models and examples of random times that are not honest times and we study the existence of classical arbitrages. This section contains two subsections that treat two different situations.

 \subsection{In a Brownian filtration: Emery's example}
We present here an example where $\tau$ is a pseudo stopping-time.

\noindent Let $S$ be defined through $dS_t=\sigma S_t dW_t$, where $W$ is a  Brownian motion and $\sigma$ a constant.
Let   $\tau = \sup \, \{t\leq
1\,:\, S_1 - 2  S_t=0\}$, that is the last time before 1 at which the
price  is equal to half of  its terminal value at time 1.

\bp In the above model   NA  holds before $\tau$.  There are classical arbitrages after $\tau$.\ep
 \proof Note that  $$ \{ \tau \leq t \}=  \{\inf_{t\leq s \leq 1} 2S_s \geq
 S_1   \} =  \{\inf_{t\leq s \leq 1} 2\frac{S_s}{S_t} \geq
 \frac{S_1}{S_t}
   \}
 $$
 Since $\frac{S_s}{S_t}, s\geq t$  and  $\frac{S_1}{S_t}$  are independent from $\F_t$, $$\P(  \inf_{t\leq s \leq 1} 2\frac{S_s}{S_t} \geq
 \frac{S_1}{S_t}
   \vert \F_t)
  = \P(  \inf_{ t\leq s \leq 1}  2 S_{s-t} \geq
   S_{1-t}) =\Phi(1-t)
  $$ where $\Phi(u)=  \P(  \inf_{  s \leq  u}  2 S_{s} \geq
   S_{u})$.
It follows that the supermartingale $Z$ is a deterministic  decreasing function,  hence, $\tau$ is a pseudo-stopping time and $S$ is a $\gg$-martingale up to time $\tau$ and there are no arbitrages up to $\tau$.

 \noindent There are obviously arbitrages after $\tau$, since, at time $\tau$, one knows the value of $S_1$ and $S_1>S_\tau$. In fact, for $t>\tau$, one has $S_t>S_\tau$, and the arbitrage occurs at any time before 1.\finproof

\subsection{In a Poisson filtration}
 This subsection develops similar examples of random times -- as in the Brownian filtration of the previous subsection -- and shows that the effects of these random times on the market's economic structure differ tremendously from the one of the previous subsection.

\noindent
In this section, we will work on a Poisson process $N$ with intensity $\lambda$ and the compensated martingale $M_t=N_t - \lambda t$. Denote
 $$T_n = \inf\{t\geq 0: N_t\geq n\}, \ \mbox{and} \ \ H^n_t= \ind_{\{T_n\leq t\}}, \  \ n=1,2.$$
 The stock price $S$ is described by
 \begin{equation}
  dS_t = S_{t-} \psi dM_t, \ \ \mbox{where}, \ \ \psi >-1, \ \mbox{and} \ \psi \neq 0.
 \end{equation}
 or equivalenty, $S_t = S_0 \exp(-\lambda \psi t + \ln(1+\psi)N_t)$.  Then,
  \begin{eqnarray*}
   M^1_t:= H^1_t - \lambda (t\wedge T_1) := H^1_t - A^1_t, \ \ \mbox{and} \ \ M^2_t:= H^2_t - \left( \lambda (t\wedge T_2) -   {\lambda (t\wedge T_1)}) \right):= H^2_t - A^2_t
  \end{eqnarray*}
are two $\mathbb{F}$-martingales.  Remark that if $\psi \in (-1,0)$, between  $T_1$ and $T_2$, the stock price increases; if $\psi >0$, between  $T_1$ and $T_2$, the stock process decreases. This would be the starting point of the existence of  arbitrages.

\subsubsection{Convex combination of two jump times}

 \noindent Below, we present an   example of random time that avoids stopping times and  the non-arbitrage property fails.

\begin{Proposition}
 Consider the random time   $\tau = k_1T_1+k_2T_2$ that  avoids $\ff$ stopping times, where $k_1+k_2 = 1$ and $k_1, k_2 >0$. Then the following properties hold:\\
 $\rm (a)$ The random time $\tau$ is not an honest time. \\
 $\rm (b) $ $\widetilde{Z}_\tau = Z_ \tau = e^{-{\lambda}k_1(T_2 - T_1)}<1,$ and $\{\widetilde{Z} =0 <Z_{-}\} =\emptyset.$\\
 $\rm (c)$ There is  a classical arbitrage before $\tau$,  given by
\begin{equation}
 \varphi_t := -e^{-\lambda \frac{k_2}{k_1}(t- T_1)}\left( \ind_{\{N_{t-} \geq 1\}} - \ind_{\{N_{t-} \geq 2 \}} \right)\frac{1}{\psi S_{t-}}\id_{\{t\leq \tau\}}.
\end{equation}
 $\rm (d)$ There exist arbitrages after $\tau$: if $\psi \in (-1,0)$, buy at $\tau$ and sell before $T_2$; if $\psi >0$, short sell at $\tau$ and buy back before $T_2.$
\end{Proposition}
\proof First, we compute the supermartingale $Z$:
  \beq \P(\tau >t \vert \F_t)&=&\ind_{T_1>t}+\ind_{\{T_1\leq t\}} \ind_{\{T_2>t\}} \P( k_1 T_1+k_2 T_2>t \vert \F_t)\eeq
  On the set  $E=( T_1\leq t )\cap (T_2>t)$, the quantity $\P( k_1 T_1+k_2 T_2>t \vert \F_t)$ is $\F_{T_1}$-measurable.
  It follows that, on $E$,
  $$
  \P( k_1 T_1+k_2 T_2>t \vert \F_t)=  \frac{\P( k_1 T_1+k_2 T_2>t , T_2>t\vert \F_{T_1})}{\P( T_2>t\vert \F_{T_1})}= \frac{e^{-\lambda\frac{1}{k_2} (t-T_1)}}{e^{-\lambda (t-T_1)}}=e^{-\lambda \frac{k_1}{k_2} (t-T_1)}, $$
  where we used the independence property of $T_1$ and $T_2 - T_1$. Therefore,
  we deduce that,
  \beq \P(\tau >t \vert \F_t)
  &=&  \ind_{\{T_1>t\}}+\ind_{\{T_1\leq t\}} \ind_{\{T_2>t\}} e^{-\lambda \frac{k_1}{k_2}  (t-T_1)}.
  \eeq

Since $Z_t= (1-H^1_t)+ H^1_t(1-H^2_t) e^{-\lambda \frac{k_1}{k_2} (t-T_1)}$, we
deduce, using the fact that $e^{-\lambda  (t-T_1)}  dH^1_t= dH^1_t$,
 \beq dZ_t&=& -dH^1_t+  e^{-\lambda \frac{k_1}{k_2}  (t-T_1)} ((1-H^2_t)dH^1_t-H^1_tdH^2_t))-  \lambda \frac{k_1}{k_2}H^1_t(1-H^2_t) e^{-\lambda \frac{k_1}{k_2}   (t-T_1)}dt\\
 &=&   e^{-\lambda   \frac{k_1}{k_2} (t-T_1)} ( -H^2 dH^1_t- H^1_t dH^2_t) - \lambda\frac{k_1}{k_2} H^1_t(1-H^2_t) e^{-\lambda  \frac{k_1}{k_2}  (t-T_1)}dt\\
 &=& - e^{-\lambda \frac{k_1}{k_2}  (t-T_1)}  dH^2_t -  \lambda\frac{k_1}{k_2} H^1_t(1-H^2_t) e^{-\lambda  \frac{k_1}{k_2}  (t-T_1)}dt\\
 &=&  dm_t - e^{-\lambda \frac{k_1}{k_2}  (t-T_1)} dA^2_t - \lambda\frac{k_1}{k_2} H^1_t(1-H^2_t) e^{-\lambda  \frac{k_1}{k_2} (t-T_1)}dt,
 \eeq
 where $$dm_t=  -  e^{-\lambda  \frac{k_1}{k_2} (t-T_1)}   dM^2_t. $$
Hence $$m_\tau= 1   -  \int_0^\tau  e^{-\lambda  \frac{k_1}{k_2}  (t-T_1)}   dM^2_t=1+   \int_{T_1}^{\tau}   e^{-\lambda \frac{k_1}{k_2}   (t-T_1)} \lambda  dt>1.$$

\noindent Now we will start proving the proposition.\\
i) Since  $\tau $  avoids stopping times,  $Z = \widetilde{Z}$.
  Note that $\widetilde{Z}_\tau = Z_ \tau = e^{-{\lambda}k_1 (T_2 - T_1)}<1$. Hence, $\tau$ is not an honest time. Since $Z>0$, we deduce that both assertions (a) and (b) hold.

 \noindent ii) Now, we will prove assertion (c). We will describe explicitly the arbitrage strategy. Note that $\{T_2 \leq t\} = \{N_t \geq 2\}$. We deduce that
 \begin{equation}
  M^2_t = \ind_{\{T_2 \leq t\}} - A^2_t =\ind_{\{N_t \geq 2\}} - A^2_t  = \ind_{\{N_{t-} \geq 1\}}\Delta N_t + \ind_{\{N_{t-} \geq 2 \}}(1 -\Delta N_t) - A^2_t.
 \end{equation}
Hence,
\begin{equation}
 \Delta M^2_t = M^2_t - \ ^p(M^2)_{t} = \left( \ind_{\{N_{t-} \geq 1\}} - \ind_{\{N_{t-} \geq 2 \}} \right)\Delta N_t  =\left( \ind_{\{N_{t-} \geq 1\}} - \ind_{\{N_{t-} \geq 2 \}} \right)\Delta M_t.
\end{equation}
Since $M^2 $ and $M $ are both purely discontinuous, we have $m_t
= 1 + (\phi \centerdot  M)_t = 1 + (\varphi \centerdot S)_t$, where
\begin{equation}
 \phi_t = -e^{-\lambda \frac{k_1}{k_2} (t- T_1)}\left( I_{\{N_{t-} \geq 1\}} - I_{\{N_{t-} \geq 2 \}} \right), \ \mbox{and }\  \varphi_t = \phi_t \frac{1}{\psi S_{t-}}.
\end{equation}

\noindent iii) { {Arbitrages after $\tau$}}: At time $\tau$, the value of $T_2$ is known for the one who has $\gg$ information. The  price process decreases before time $T_2$, however, waiting up time $T_2$ does not lead to an arbitrage
 Setting $\Delta =T_2-\tau$ (which is known at time $\tau$), there is an arbitrage selling short $S$ at time $\tau$ for a delivery at time $\tau+\frac 12 \Delta$.  The strategy is admissible, since between $T_1$ and $T_2$,   the quantity $S_t$ is bounded by $S_0 (1+\varphi)$. This ends the proof of the proposition.\finproof

\subsubsection{Minimum of two scaled jump times}
We give now an example of a  non honest random time, which
does  not avoid $\ff$ stopping time and induces classical
arbitrage opportunities.

 \begin{Proposition}
  Consider the same market as before, and define  $\tau= T_1 \wedge aT_2$, where $0<a<1$. Then, the following properties hold:\\
  $\rm (a)$ $\tau$ is not an honest time and  does not avoid $\mathbb{F}$-stopping times,\\
  $\rm (b)$ $ Z_\tau = \ind_{\{T_1 > aT_2\}}e^{-\beta a T_2} (\beta a T_2 +1)<1$ and $\widetilde{Z}_\tau = e^{-\beta a T_2} (\beta a T_2 +1)<1$,  and $\{\widetilde{Z} =0 <Z_{-}\} =\emptyset.$\\
  $\rm (c)$ There exists a classical arbitrage before $\tau$  given by
  \begin{equation}
 \varphi_t = -e^{-\beta t}(\beta t+1) \left( \ind_{\{N_{t-} \geq 0\}} -\ind_{\{N_{t-} \geq 1 \}} \right) \frac{1}{\psi S_{t-}}\id_{\{t\leq \tau\}}, \ \ \mbox{where} \ \, \beta= \lambda (1/a -1).
\end{equation}
$\rm (d)$ There exist arbitrages after $\tau$: if $\psi \in (-1,0)$, buy at $\tau$ and sell before $\tau/a$; if $\psi >0$, short sell at $\tau$ and buy back before $\tau/a.$
 \end{Proposition}
\proof First, let us compute the supermartingale $Z$,
 \beq Z_t&=&\ind_{\{T_1>t\}}\P( aT_2>   t  \vert \F_t)=\ind_{T_1>t}\frac{\P( aT_2>   t , T_1>t)}{ \P(T_1>t)} \\
 &=&  \ind_{\{T_1>t\}} e^{\lambda t}\,\E( \ind _{T_1>t}  e^{-\lambda (\frac ta -T_1)^+})\\
 &=& \ind_{\{T_1>t\}} e^{\lambda t}\int_t^{t/a}   e^{-\lambda (\frac ta -x)}\lambda e^{-\lambda x}dx + \ind_{\{T_1>t\}} e^{\lambda t}\int_{t/a}^{\infty} \lambda e^{-\lambda y}dy\\
 &=& \ind_{\{T_1>t\}} e^{- \beta t} (\beta t + 1),
 \eeq
 where $\beta= \lambda (1/a -1)$. In particular $Z_\tau = \ind_{T_1 > aT_2}e^{-\beta a T_2} (\beta a T_2 +1)<1$.  Similar computation as above leads to $\widetilde Z_t=Z_{t-} =\ind_{\{T_1\geq t\}} e^{- \beta t} (\beta t + 1)$. This proves assertions (a) and (b). \\
 i) Here, we will prove assertion (c). Thanks to It\^o's formula, we have
 \begin{eqnarray*}
  dZ_t =  -e^{-\beta t}(\beta t+1) dH^1_t - \ind_{t\leq T_1}  \beta^2 e^{-\beta t} tdt = -e^{-\beta t}(\beta t +1) dM^1_t  - e^{-\beta t}(\beta t +1) dA^1_t  - \ind_{t\leq T_1}  \beta^2 e^{-\beta t} tdt.
 \end{eqnarray*}
Therefore,
 $$
 dm_t=-e^{-\beta t}(\beta t  +1) dM^1_t.
 $$
Hence \begin{eqnarray}
m_\tau&=&1+\ind_{\{aT_2<T_1\}} \lambda \left ( \frac{2(1-e^{-\beta aT_2})}{\beta}-aT_2 e^{-\beta aT_2}\right) \nonumber \\
&&+  \ind_{\{T_1<aT_2\}} \left ( 2\lambda \frac{1-e^{-\beta  T_1}}{\beta}- \lambda T_1 e^{-\beta  T_1}-T_1 \beta e^{-\beta  T_1} - e^{-\beta T_1} +\ind_{\{T_1 =0\}}\right)
\end{eqnarray} and, using the fact that when $x>0$,  $1-e^{-x}-xe^{-x}>0$  and $2\lambda \beta e^{-\beta x} - \lambda \beta xe^{-\beta x} - x\beta^2 e^{-\beta x} - \beta e^{-\beta x} + \beta I_{x=0} >0$, one obtains $m_\tau >1$;
hence the existence of classical  arbitrages.

\noindent Now, we describe explicitly the arbitrage strategy. Notice that  $\{T_1 \leq t\} = \{N_t \geq 1\}$. We deduce that
 \begin{equation}
  M^1_t = \ind_{\{T_1 \leq t\}} - A^1_t = \ind_{\{N_t \geq 1\}} - A^1_t  =\ind_{\{N_{t-} \geq 0\}}\Delta N_t + \ind_{\{N_{t-} \geq 1 \}}(1 -\Delta N_t) - A^1_t.
 \end{equation}
Hence,
\begin{equation}
 \Delta M^1_t = M^1_t - \ ^p(M^1)_{t} = \left( \ind_{\{N_{t-} \geq 0\}} - \ind_{\{N_{t-} \geq 1 \}} \right)\Delta N_t  =\left( I_{\{N_{t-} \geq 0\}} - I_{\{N_{t-} \geq 1 \}} \right)\Delta M_t.
\end{equation}
Since $M^1 $ and $M $ are both purely discontinuous, we have
$m  = 1 + \varphi \centerdot S$, where
\begin{equation}
 \varphi_t = -e^{-\beta t}(\beta t +1) \left( \ind_{\{N_{t-} \geq 0\}} - \ind_{\{N_{t-} \geq 1 \}} \right)\frac{1}{\psi S_{t-}}.
\end{equation}
ii) The proof of assertion (d) follows the same proof of assertion (d) of Proposition 5.6. This ends the proof of the proposition.
\finproof

\subsubsection{Maximum of two scaled jump times}

 \begin{Proposition}
  Consider the same market as before. and define $\tau= T_1 \vee aT_2$, where $0<a<1$. Then, the following properties hold:\\
  $\rm (a)$ $\tau$ is not an honest time and  does not avoid $\mathbb{F}$-stopping times. \\
  $\rm (b)$ $Z_\tau = \frac{\lambda \tau e^{-\lambda {\tau}/{a}}}{1 - e^{-\lambda \tau }} <1,$ and $\widetilde{Z}_\tau = I_{\{T_1 \geq aT_2\}} + I_{\{T_1 <  aT_2\}}Z_\tau \not\equiv 1,$  and $\{\widetilde{Z} =0 <Z_{-}\} =\emptyset.$\\
  $\rm (c)$ There exists  a classical arbitrage before $\tau$   given by
  \begin{equation}
 \varphi_t = - \left( 1 - \frac{\lambda t e^{-\lambda \frac{t}{a}}}{1 - e^{-\lambda t}} \right)\left( I_{\{N_{t-} \geq 0\}} - I_{\{N_{t-} \geq 1 \}} \right) \frac{1}{\psi S_{t-}}.
\end{equation}
 $\rm (d)$ There exist classical arbitrages after $\tau$: if $\psi \in (0,1)$ and $T_1<aT_2$, buy at $\tau$ and sell before $\tau/a$; if $\psi >0$ and $T_1<aT_2$, short sell at $\tau$ and buy back before $\tau/a.$
 \end{Proposition}
\proof First, let us compute the supermartingale $Z$,
\begin{eqnarray*}
 1-Z_t &=& P\left( \tau \leq t | {\cal F}_t \right)  = P\left( T_1\vee a T_2\leq t | {\cal F}_t \right) \\
 &=& \ind_{\{T_1 \leq t\}} P\left(T_2\leq \frac{t}{a} \Big | {\cal F}_t \right) = \ind_{\{T_1 \leq t\}} \frac{P\left(T_2\leq \frac{t}{a}, T_1 \leq t\right)}{P\left( T_1 \leq t \right)} \\
 &=& \ind_{\{T_1 \leq t\}} \frac{1}{1 - e^{-\lambda t}} \int_0^t\left( 1 - e^{-\lambda \left(  \frac{t}{a} -y \right)} \right)\lambda e^{-\lambda y}dy \\
 &=& \ind_{\{T_1 \leq t\}}\left( 1 - \frac{\lambda t e^{-\lambda \frac{t}{a}}}{1 - e^{-\lambda t}} \right).
 \end{eqnarray*}
 Therefore
 \begin{equation}
  Z_t = 1 -  \ind_{\{T_1 \leq t\}}\left( 1 - \frac{\lambda t e^{-\lambda \frac{t}{a}}}{1 - e^{-\lambda t}} \right) = \ind_{\{T_1 > t\}} + I_{\{T_1 \leq t\}}\frac{\lambda t e^{-\lambda \frac{t}{a}}}{1 - e^{-\lambda t}},
 \end{equation}
and $$Z_\tau = \frac{\lambda (T_1\vee a T_2) e^{-\lambda \frac{T_1\vee a T_2}{a}}}{1 - e^{-\lambda (T_1\vee a T_2)}} <1.$$
 Using the same type of arguments give
\begin{equation}
  \widetilde{Z}_t =  \ind_{\{T_1 \geq  t\}} + \ind_{\{T_1 < t\}}\frac{\lambda t e^{-\lambda \frac{t}{a}}}{1 - e^{-\lambda t}}, \ \ \mbox{and} \ \ \widetilde{Z}_\tau = \ind_{\{T_1 \geq aT_2\}} + \ind_{\{T_1 <  aT_2\}} \frac{\lambda (T_1\vee a T_2) e^{-\lambda \frac{T_1\vee a T_2}{a}}}{1 - e^{-\lambda (T_1\vee a T_2)}}.
 \end{equation}
 These allows us to conclude that both assertions (a0 and (b) hold. The proof of the assertion (d) is similar to that asertion (d) of the previous proposition. The remaining part of the proof will address assertion (c).\\
 By putting $K_t = 1 - \frac{\lambda t e^{-\lambda \frac{t}{a}}}{1 - e^{-\lambda t}} $, $Z_t=1-H^1_tK_t$ and applying It\^o formula, we derive that
 \begin{eqnarray}
  dZ_t = - d(H^1K)_t = -K_t dH^1_t - H^1_tdK_t = -K_t dM^1_t  - K_t dA^1_t - H^1_tdK_t.
 \end{eqnarray}
Hence,
\begin{equation}
 m_t = 1 - \int_0^t K_t dM^1_t.
\end{equation}

\noindent Now, we describe explicitly the arbitrage strategy. Notice that  $\{T_1 \leq t\} = \{N_t \geq 1\}$. We deduce that
 \begin{equation}
  M^1_t = I_{\{T_1 \leq t\}} - A^1_t = I_{\{N_t \geq 1\}} - A^1_t  = I_{\{N_{t-} \geq 0\}}\Delta N_t + I_{\{N_{t-} \geq 1 \}}(1 -\Delta N_t) - A^1_t.
 \end{equation}
Hence,
\begin{equation}
 \Delta M^1_t = M^1_t - \ ^p(M^1)_{t} = \left( I_{\{N_{t-} \geq 0\}} - I_{\{N_{t-} \geq 1 \}} \right)\Delta N_t  =\left( I_{\{N_{t-} \geq 0\}} - I_{\{N_{t-} \geq 1 \}} \right)\Delta M_t.
\end{equation}
Since $M^1_t$ and $M_t$ are both purely discontinuous, we have
$m_t = 1 + \varphi \cdot S_t$, where
\begin{equation}
 \varphi_t = - K_t \left( I_{\{N_{t-} \geq 0\}} - I_{\{N_{t-} \geq 1 \}} \right)\frac{1}{\psi S_{t-}}.
\end{equation}
\finproof

\section{ NUPBR for particular models}
In this section, we address some interesting practical models, for which we prove that the NUPBR remains valid up to $\tau$. The originality of this part -- as we mentioned in the introduction and the abstract -- lies in the simplicity of the proof. A general and complete analysis about the NUPBR is addressed in full generality in Choulli et al. (2013). Throughout this section, we will assume that $Z>0$.

\subsection{Before $\tau$}

Let $\widehat m$ be the $\gg$-martingale  stopped at time $\tau$
associated with $m$ by (\ref{jeulinavant}), on $\{t \leq \tau\}$
$$\widehat{m}_t := m_t ^\tau  - \int_{0}^{t } { \frac{d \langle m,
m\rangle^{\mathbb F}_s }{Z_{s }} }\,.$$

\subsubsection{Case of continuous filtration}
We start with the particular case of
continuous martingales and prove that, for any random time $\tau$,
  NUPBR holds  before $\tau$.

\noindent  We note that the continuity assumption implies that the martingale part of $Z$ is continuous and   that   the optional and Doob-Meyer decompositions of $Z$ are the same.
 \bp Assume that all $\ff$-martingales are continuous.
Then, for any random time $\tau$,
NUPBR holds before $\tau$. A $\gg$-local martingale deflator for $S^\tau$ is given by  $dL_t= -\frac{L_t}{Z_{t }}d\widehat m_t.$  \ep
 \proof  We make a use of Theorem \ref{deflator} and we provide a $\gg$-local martingale deflator for $S^\tau$. Define the positive $\gg$-local martingale $L$ as $dL_t= -\frac{L_t}{Z_{t }}d\widehat m_t.$
 Then, if  $SL$ is a $\gg$-local martingale, NUPBR holds. Recall that, using
(\ref{jeulinavant}) again,
$$\widehat{S}_t  := S_t^\tau  - \int_{0}^{t  \wedge \tau} { \frac{d \langle S,
m\rangle^{\mathbb F}_s }{Z_{s }} }$$ is a $\gg$-local martingale.
 From integration by parts, we obtain (using  that the bracket of
continuous martingales does not depend on the filtration)
\beq d(LS^{\tau})_t&=& L_{t }dS^{\tau}_t+S_{t }dL_t+d\cro{L,S^{\tau}}^\gg_t\\
 &\mart&  L _{t }   \frac{1}{Z_t}d \cro{S,m}^\ff_t+   \frac{1}{Z_{t}} L_{t}  d\cro{S,\widehat m}^\gg_t \\
& \mart &    L_{t }    \frac{1}{Z_t} \left ( d \cro{S,m}_t- d\cro{S,
m }_t \right)=0\eeq where $X\mart Y$  is a notation for $X-Y$ is a $\gg$-local martingale. \finproof

\brem
\label{s1}
If $\tau$ is an honest time and Predictable Representation Property holds with respect to $S$
 then, as a consequence of Theorem \ref{sa}, the NA condition does not hold, hence NFLVR condition does not hold neither. That in turn implies that all the $\gg$-local martingale deflators for $S^\tau$ are strict $\gg$-local martingales.
\erem

\subsubsection{Case of a Poisson filtration}

 We assume that $S$ is an $\ff$-martingale of the form $dS_t=S_{t-}\psi_t dM_t$, where  $\psi$ is a
predictable process, satisfying $\psi>-1$ and $\psi \neq 0$, where
$M$ is the compensated martingale of a standard Poisson process.

\noindent In a Poisson setting, from Predictable Representation Property,  $dm_t= \nu_tdM_t$ for some $\ff$-predictable process $\nu$,  so that, on $t\leq \tau$, $$d \widehat m_t= dm_t -\frac {1}{Z_{t-}}d\cro{m,m}_t=dm_t -\frac {1}{Z_{t-}} \lambda \nu^2_tdt$$

\bp \label{poissoncase}In a Poisson setting, for any random time $\tau$,  NUPBR holds
before $\tau$ since
$$L= {\cal E}\left(-\frac{1}{Z_{-}+\nu}\centerdot \widehat m\right)= {\cal E}\left(-\frac{\nu}{Z_{-}+\nu}\centerdot\widehat M\right),$$     is a $\gg$-local martingale deflator for $S^{\tau}$.\ep
\proof  We make a use of Theorem \ref{deflator} and we are looking for a  $\gg$-local martingale deflator of the form $dL_t=L_{t-}\kappa_t d\widehat m_t$ (and $\psi_t \kappa_t>-1$) so that $L$ is positive and $S^\tau L$ is a $\gg$-local martingale. Integration by parts formula leads to (on $t\leq \tau$)
  \beq d(LS)_t&=&L_{t-}dS_t+S_{t-}dL_t+d[L,S]_t \\&\mart& L_{t-} S_{t-}\psi_t \frac{1}{Z_{t-}} d\cro {M,m}_t + L_{t-}S_{t-} \kappa_t \psi_t \nu_t  dN_t\\
&\mart& L_{t-} S_{t-}\psi_t \frac{1}{Z_{t-}} \nu_t \lambda dt + L_{t-}S_{t-} \kappa_t \psi_t \nu_t   \lambda (1+ \frac 1{Z_{t-}} \nu_t)dt \\&=&  L_{t-} S_{t-}\psi_t \nu_t   \lambda \left( \frac{1}{Z_{t-}}     +   \kappa_t   (1+ \frac 1{Z_{t-}} \nu_t)
\right)dt .\eeq
Therefore, for $\kappa_t =-\frac{1}{Z_{t-}+\nu_t}$, one obtains a
deflator. Note that $$dL_t=L_{t-}\kappa_t d\widehat
m_t=-L_{t-}\frac{1}{Z_{t-}+\nu_t}\nu_td\wh M_t$$  is indeed a
positive $\gg$-local martingale, since   {$\frac{1}{Z_{t-}+\nu_t}\nu_t<1$}.
\finproof

\brem
\label{s3}
If $\tau$ is an honest time and Predictable Representation Property holds with respect to $S$
 then all the $\gg$-local martingale deflators for $S^\tau$ are strict $\gg$-local martingales.
\erem

\subsubsection{L\'evy processes}

Assume that $S = \psi \star (\mu -\nu)$ where $\mu$ is   the jump measure of a L\'evy process  and $\nu$ its compensator. Here, $\psi \star (\mu -\nu)$ is the process
 $\int_0^\cdot \int \psi(x,s) (\mu(dx,ds)-\nu (dx,ds))$.
The martingale $m$ admits a representation as $m= \psi^m \star (\mu-\nu) $.
Then, using \eqref{jeulinavant}, the $\gg$-compensator of $\mu$ is $\nu ^\gg$ where
$$\nu^{\mathbb G}(dt,dx)= \frac{1}{Z_{t-}}\left(Z_{t-}+\psi^m(t,x)\right)  \nu(dt,dx) $$ i.e.,  $S$ admits a $\gg$-semimartingale decomposition of the form
  $$S = \psi \star (\mu -\nu^\gg) -\psi \star (\nu-\nu ^\gg)$$

\bp Consider the positive $\gg$-local martingale
$$
L:={\cal E}\left(- \frac{\psi^m}{Z_-+\psi^m}I_{\Lbrack 0,\tau\Lbrack}\star(\mu-\nu^{\gg})\right).$$
Then $L$ is a $\gg$-local martingale deflator for $S^{\tau}$, and hence $S^{\tau}$ satisfies NUPBR.
\ep
\proof  We make a use of Theorem \ref{deflator} and our goal is to find a  positive $\gg$-local martingale $L$ of the
form
$$dL_t=L_{t-}\kappa_t d\widehat m_t$$ so that $LS^\tau$ is a
$\gg$-local martingale.

\noindent  From integration by parts formula

\beq d(SL) &\mart&- L_- \psi\star (\nu-\nu ^\gg)+ d[S,L] = - L_- \psi\star (\nu-\nu ^\gg)+ L_ - \psi  \psi^m \kappa \star \mu\\
&\mart& - L_-  \psi \star (\nu-\nu ^\gg)+ L_- \psi  \psi^m \kappa  \star \nu ^\gg\\
&=&- L_- \psi \left( 1- (1+\psi^m \kappa ) \frac{1}{Z_{ -}}\left(Z_{ -}+\psi^m \right)  \right)  \star \nu \eeq
Hence the possible choice $\kappa =- \frac{1}{Z_-+\psi^m}$.
 {It can be checked that indeed, $L$ is a positive $\gg$-local martingale See \cite{acdj}.}
\finproof

\subsection{After $\tau$}
We now assume that $\tau$ is an honest time, which satisfies
$Z_\tau <1$ (for integrability reasons). This condition and Lemma \ref{stricthonest} imply in particular that $\tau$ does not avoids $\ff$-stopping times. For the further discussion on the condition $Z_\tau<1$ we refer the reader to \cite{thin}. Note also  that,
in the case of continuous filtration, and $Z_\tau=1$, NUPBR fails
to hold after $\tau$ (see \cite{fjs}).

\noindent After \eqref{after}, for any $\ff$-martingale $X$  (in particular for $m$ and $S$)
$$\widehat{X}_t  := X_t^\tau  - \int_{0}^{t  \wedge \tau} { \frac{d \langle X,
m\rangle^{\mathbb F}_s }{Z_{s }} } + \int_{t  \wedge \tau }^t{} {
\frac{d \langle X, m\rangle^{\mathbb F}_s }{1- Z_{s }} }$$ is a
$\gg$-local martingale.

\subsubsection{Case of continuous filtration}
We start with the particular case of
continuous martingales and prove that, for any honest time $\tau$ such that $Z_\tau<1$, NUPBR holds  after $\tau$.

 \bp Assume that  $\tau$ is an honest time, which satisfies $Z_\tau <1$  and that   all $\ff$-martingales are continuous. Then, for any honest  time $\tau$,
NUPBR  holds  after $\tau$.   A $\gg$-local martingale deflator for $S-S^\tau$ is given by $dL_t= -\frac{L_t}{1- Z_{t }}d\widehat m_t.$\ep
 \proof  We use Theorem \ref{deflator} as usual. The proof is based on It\^o's calculus. Looking for a $\gg$-local martingale deflator of the form
$dL_t=L_t\kappa_t d\wh m_t$, and using
  integration by parts formula, we obtain  that, for $\kappa= - (1- Z  )^{-1}$, the process  $L(S-S^\tau)$ is a $\gg$-local martingale. \finproof

\brem
\label{s2}
If Predictable Representation Property holds with respect to $S$
 then, as a consequence of Theorem \ref{sa}, the NA condition does not hold, hence NFLVR condition does not hold neither. That in turn implies that all the $\gg$-local martingale deflators for $S-S^\tau$ are strict $\gg$-local martingales.
\erem

\subsubsection{Case of a Poisson filtration}
We assume that $S$ is an $\ff$-martingale  of the form
$dS_t=S_{t-}\psi_t dM_t$, with  $\psi$ is a predictable process,
satisfying $\psi>-1$.

\noindent The decomposition formula \eqref{after} reads after $\tau$ as

 $$ \widehat S_t= (\id_{]\tau,\infty[}\cdot S)_t +\int_{t\lor \tau}^t \frac {1}{1-Z_{s-}}d\cro{S,m}_s=(\id_{]\tau,\infty[}\cdot S)_t +\lambda\int_{t\lor \tau}^t \frac {1}{1-Z_{s-}} \nu_s \psi_s S_{s-} ds.$$

\bp \label{poissoncase} Let $\ff$ be a Poisson filtration and
$\tau$ is an honest time satisfying $Z_\tau<1$. Then,  NUPBR holds
after $\tau$ since $$L= {\cal
E}\left(\frac{1}{1-Z_{-}-\nu}\centerdot \widehat m\right)= {\cal
E}\left(\frac{\nu}{1-Z_{-}-\nu}\id_{]\tau,\infty[}\centerdot\widehat
M\right),$$ is a $\gg$-local martingale deflator for $S-S^{\tau}$.
 \ep
\proof  We make a use of Theorem \ref{deflator} and we are looking for a $\gg$-local martingale deflator of the form
$dL_t=L_{t-}\kappa_t d\widehat m_t$ (and $\psi_t \kappa_t>-1$) so
that $L$ is positive $\gg$-local martingale and $(S-S^\tau) L$ is
a $\gg$-local martingale. Integration by parts formula leads to
  \beq d(L(S-S^\tau))_t&=&L_{t-}d(S-S^\tau)_t+(S_{t-}-S^\tau_{t-})dL_t+d[L,S-S^\tau]_t \\&\mart& -\lambda L_{t-} S_{t-}\nu_t \psi_t \frac{1}{1-Z_{t-}} \id_{\{t>\tau\}}dt + L_{t-}S_{t-} \kappa_t \psi_t \nu_t  \id_{\{t>\tau\}} dN_t\\
&\mart& -\lambda L_{t-} S_{t-}\nu_t \psi_t \frac{1}{1-Z_{t-}} \id_{\{t>\tau\}}dt + \lambda L_{t-}S_{t-} \kappa_t \psi_t \nu_t  \id_{\{t>\tau\}} (1- \frac 1{1-Z_{t-}} \nu_t)dt \\&=& \lambda L_{t-}S_{t-} \psi_t \nu_t  \id_{\{t>\tau\}}\left (-\frac 1{1-Z_{t-}}  +   \kappa_t   (1- \frac 1{1-Z_{t-}} \nu_t) \right )dt .\eeq

\noindent Therefore, for $\kappa_t =\frac{1}{1-Z_{t-}-\nu_t}$, one obtains a $\gg$-local martingale
deflator. Note that $$dL_t=L_{t-}\kappa_t d\widehat
m_t=L_{t-}\frac{1}{1-Z_{t-}-\nu_t}\nu_t\id_{\{t>\tau\}}d\wh M_t$$  is indeed a
positive $\gg$-local martingale, since   {$\frac{1}{1-Z_{t-}-\nu_t}\nu_t\Delta N_t>-1$}.
\finproof

\brem
\label{s4}
If Predictable Representation Property holds with respect to $S$
 then, all the $\gg$-local martingale deflators for $S-S^\tau$ are strict $\gg$-local martingales.
\erem

\subsubsection{L\'evy processes}

Assume that $S = \psi \star (\mu -\nu)$ where $\mu$ is   the jump measure of a L\'evy process  and $\nu$ its $\ff$-compensator.
Then, by \eqref{after}, the $\gg$-compensator of $\mu$ is $\nu ^\gg$ where
$$\nu^{\mathbb G}(dt,dx)= \left(1+\id_{\{t\leq \tau\}}\frac{1}{Z_{t-}}\psi^m(t,x)-\id_{\{t>\tau\}}\frac{1}{1-Z_{t-}}\psi^m(t,x)\right)  \nu(dt,dx) $$ i.e.,  $S$ admits a $\gg$-semimartingale decomposition of the form
  $$S = \psi \star (\mu -\nu^\gg) -\psi \star (\nu-\nu ^\gg)$$

\bp Assume that $\tau$ be an honest time satisfying $Z_\tau<1$  in
a L\'evy framework. Then, the positive $\gg$-local martingale
$$
L:={\cal E}\left( \frac{\psi^m}{1-Z_--\psi^m}I_{\Lbrack
\tau,\infty\Rbrack}\star(\mu-\nu^{\gg})\right),$$
  is a $\gg$-local martingale deflator for $S-S^{\tau}$, and hence $S-S^{\tau}$ satisfies NUPBR.
\ep
\proof  We use of Theorem \ref{deflator} again. Our goal is to find a  positive $\gg$-local martingale
$L$ of the form $$dL_t=L_{t-}\kappa_t \id_{\{t>\tau\}}d\widehat
m_t$$ so that $L(S-S^\tau)$ is a $\gg$-local martingale.

\noindent  From integration by parts formula
\beq d(L(S-S^\tau)) &\mart&- L_-d (S-S^\tau)+ d[S,L] \\
&=& - L_- \psi\frac{\psi^m}{1-Z_-}\id_{]\tau,\infty[}\star \nu+ L_- \kappa \psi  \psi^m \id_{]\tau,\infty[}\star \mu\\
&\mart& - L_- \psi\frac{\psi^m}{1-Z_-}\id_{]\tau,\infty[}\star \nu+ L_- \kappa \psi  \psi^m \id_{]\tau,\infty[}\star \nu ^\gg\\
&=&- L_- \psi \psi^m
\id_{]\tau,\infty[}\left(-\frac{1}{1-Z_-}+\kappa
(1-\frac{\psi^m}{1-Z_-})  \right)  \star \nu \eeq Hence the
possible choice $\kappa =\frac{1}{1-Z_--\psi^m}$. \finproof

\section*{Conclusions}

In this paper we have treated the question whether the no-arbitrage conditions are stable with respect to progressive enlargement of filtration. We focused on two components of No Free Lunch with Vanishing Risk concept, namely on No Arbitrage Opportunity and No Unbounded Profit with Bounded Risk. The problem was divided into stability before and after random time containing extra information.

\noindent The question regarding No Arbitrage Opportunity condition was answered in the case of Brownian filtration and Poisson filtration for special case of honest time, moreover particular examples of non-honest times were described.
Both, Brownian and Poisson filtrations possess an important, and crucial from our problem point of view, characteristic of Predictable Representation Property. One may further investigate similar problem without assuming market completness. One may as well consider other example/classes of non-honest random times.

\noindent Afterwards, we handeled with stability of NUPBR concept in very particular situations, namely in continuous martingale case, standard Poisson process case and L\'evy process case. We provided reults with simple proofs in those particular situations.
We emphasize again that in full generality the problem is solved in \cite{acdj} revealing as well results within progressive enlargement of filtration theory.

\noindent Combining results on NA and NUPBR conditions we concluded (in Remarks \ref{s1}, \ref{s3}, \ref{s2}, \ref{s4}) that some $\gg$-local martingales are in fact $\gg$-strict local martingales. That provides a way to construct strict local martingale in enlarged Brownian and Poisson filtrations.

\appendix
\section{Appendix}

 Let $(A_t,\,t\geq 0)$ be an integrable increasing process (not
necessarily $\ff$-adapted). There exists a unique integrable
$\ff$-optional increasing process $(A^{o }_t,t\geq 0)$, called the
dual optional  projection of $A$ such that
$$\E\left( \int_{[0,\infty[} U_sdA_s \right)=\E\left( \int_{[0,\infty[}
U_sdA_s ^{o }\right)$$ for any positive $\ff$-optional process
$U$.\\ There exists a unique integrable $\ff$-predictable
increasing process $(A^{p }_t,t\geq 0)$, called the dual
predictable projection of $A$ such that
$$\E\left( \int_{[0,\infty[} U_sdA_s \right)=\E\left( \int_{[0,\infty[}
U_sdA_s ^{p }\right)$$ for any positive $\ff$-predictable process
$U$.
\\

\noindent Acknowledgement: { This research benefited from the support of the "Chair Markets in Transition", under the aegis
of Louis Bachelier laboratory, a joint initiative of Ecole Polytechnique, Universit\'e d'Evry Val
d'Essonne and Fed\'eration Bancaire Fran\c{c}aise.}

\end{document}